\documentclass[12pt,epsfig,epsf]{amsart}
\usepackage{amssymb}
\usepackage{url}
\usepackage{epic, eepic, subfigure, floatflt}
\usepackage{amsfonts, epsfig, latexsym, amsmath}
\usepackage{color}
\def\QuotS#1#2{\leavevmode\kern-.0em\raise.2ex\hbox{$#1$}\kern-.1em/\kern-.1em\lower.25ex\hbox{$#2$}}
\DeclareMathOperator{\tdeg}{deg}

\title{Combinatorial cube packings in cube and torus}

\author{Mathieu Dutour Sikiri\'c}
\address{Mathieu Dutour Sikiri\'c, Rudjer Bo\u skovi\'c Institute, Bijenicka 54, 10000 Zagreb, Croatia}
\email{mdsikir@irb.hr}

\thanks{The first author was partly supported by the Croatian Ministry of Science, Education and Sport under contract 098-0982705-2707}

\author{Yoshiaki Itoh}
\address{The Institute of Statistical Mathematics, 4-6-7 Minami-Azabu, Minato-ku, Tokyo 106-8569, Japan}
\email{itoh@ism.ac.jp}

\DeclareMathOperator{\Sym}{Sym}
\DeclareMathOperator{\tdim}{dim}

\DeclareMathOperator{\tord}{ord}
\newcommand{\RR}{\ensuremath{\mathbb{R}}}

\newcommand{\QQ}{\ensuremath{\mathbb{Q}}}

\newcommand{\ZZ}{\ensuremath{\mathbb{Z}}}

\newtheorem{theorem}{Theorem}[section]

\newtheorem{conjecture}[theorem]{Conjecture}

\newtheorem{lemma}[theorem]{Lemma}

\newtheorem{proposition}[theorem]{Proposition}

\begin{document}
\maketitle

\begin{abstract}
We consider sequential random packing
of cubes $z+[0,1]^n$ with $z\in \frac{1}{N}\ZZ^n$ into the cube $[0,2]^n$
and the torus $\QuotS{\RR^n}{2\ZZ^n}$ as $N\to\infty$.
In the cube case $[0,2]^n$ as $N\to\infty$ the random cube packings
thus obtained are reduced
to a single cube with probability $1-O\left(\frac{1}{N}\right)$.
In the torus case the situation is different: for $n\leq 2$,
sequential random cube packing yields cube tilings, but for $n\geq 3$
with strictly positive probability, one obtains non-extensible cube packings.

So, we introduce the notion of combinatorial cube packing, which instead
of depending on $N$ depend on some parameters.
We use use them to derive an expansion of the packing density in 
powers of $\frac{1}{N}$. The explicit computation is done in the cube
case.
In the torus case, the situation is more complicate and 
we restrict ourselves to the case $N\to\infty$ of strictly positive
probability.
We prove the following results for torus combinatorial cube packings:
\begin{itemize}
\item We give a general Cartesian product construction.
\item We prove that the number of parameters is at least $\frac{n(n+1)}{2}$ and we conjecture it to be at most $2^n-1$.
\item We prove that cube packings with at least $2^n-3$ cubes are extensible.
\item We find the minimal number of cubes in non-extensible cube packings
for $n$ odd and $n\leq 6$.
\end{itemize}

\end{abstract}

\section{Introduction}
Two cubes $z+[0,1]^n$ and $z'+[0,1]^n$ are 
{\em non-overlapping} if the relative interiors
$z+]0,1[^n$ and $z'+]0,1[^n$ are disjoints.
A family of cubes $(z^i+[0,1]^n)_{1\leq i\leq m}$ with
$z^i\in \frac{1}{N} \ZZ^n$ and
$N\in \ZZ_{>0}$ is called a {\em discrete cube packing} if any two cubes
are non-overlapping.
We consider packing of cubes $z+[0,1]^n$ with $z\in \frac{1}{N}\ZZ^n$ into the cube $[0,2]^n$ and the torus $\QuotS{\RR^n}{2\ZZ^n}$.
In those two cases, two cubes $z+[0,1]^n$ and $z'+[0,1]^n$ are non-overlapping if and only if there exist an index $i\in \{1,\dots,n\}$ such that $z_i\equiv z'_i+1\pmod 2$.
A discrete cube packing is a {\em tiling} if the number of cubes is $2^n$
and it is {\em non-extensible} if it is maximal by inclusion with less than
$2^n$ cubes.

A {\em sequential random cube packing} consists of putting a cube $z+[0,1]^n$
with $z\in \frac{1}{N}\ZZ^n$ uniformly at random in
the cube $[0,2]^n$ or the torus $\QuotS{\RR^n}{2\ZZ^n}$
until a maximal packing is obtained.
Let us denote by $M_N^C(n)$, $M_N^T(n)$ the random variables of number
of cubes of those non-extensible cube packings and by $E(M_N^C(n))$,
$E(M_N^T(n))$ their expectation.
We are interested in the limit $N\to\infty$ and we prove that if $N>1$ then
\begin{equation}\label{ExpansionCombCubPack}
E(M_N^U(n))=\sum_{k=0}^{\infty} \frac{U_{k}(n)}{(N-1)^k}
\mbox{~with~}U\in \{C, T\}\mbox{~and~}U_k(n)\in \QQ
\end{equation}
In the cube case we prove that $C_k(n)$ are polynomials of degree $k$, which
we compute for $k\leq 6$ (see Theorem \ref{MainTheoremRigidBoundaryCase}).
In particular, $C_{0}=1$, since as $N\to\infty$ with
probability $1-O(\frac{1}{N})$, one cannot add any more cube after the
first one.
In the torus case the coefficients $T_{k}(n)$ are no longer polynomials in
the dimension $n$.
The first coefficient $T_{0}(n)=\lim_{N\to\infty} E(M_N^T(n))$ is known
only for $n\leq 4$ (see Table \ref{DataInformation}).
But we prove in Theorem \ref{UpperBoundLamination}
that if $n\geq 3$ then $T_0(n) < 2^n$. This upper bound is related to
the existence in dimension $n\geq 3$ of non-extensible torus cube packings
(see Figure \ref{EnumerationDim3}, Table \ref{DataInformation}, 
Theorem \ref{UpperBoundLamination} and Section \ref{NonExtensibilityQuestions}).

Those results are derived using the notion of {\em combinatorial cube
packings} which is introduced in Section \ref{SecCombCubPacking}.
A combinatorial cube packing does not depend on $N$ but instead on
some parameters $t_i$; to a cube or torus discrete cube
packing ${\mathcal{CP}}$,
one can associate a combinatorial cube packing
${\mathcal{CP}}'=\phi({\mathcal{CP}})$.
Given a combinatorial cube packing ${\mathcal{CP}}$ the probability
$p({\mathcal{CP}}, N)$ of obtaining a discrete cube packing ${\mathcal{CP}}'$
with $\phi({\mathcal{CP}}')={\mathcal{CP}}$ is a fractional function of $N$.
We say that ${\mathcal{CP}}$ is obtained with {\em strictly positive
probability} if the limit $\lim_{N\to\infty}p({\mathcal{CP}}, N)$
is strictly positive.

In Section \ref{DiscreteCubePackCube} the method of combinatorial cube
packings is applied to the cube case and the polynomials $C_k$ are
computed for $k\leq 6$. In the torus case, the situation is more complicated
and we restrict ourselves to the case of strictly positive probability, 
i.e. the limit case $N\to\infty$.
In Section \ref{LaminationTorusCubePacking} we consider a Cartesian product
construction for continuous cube packings obtained with strictly positive
probability. The related lamination construction is used to derive an upper
bound on $E(M_{\infty}^T(n))$ in Theorem \ref{UpperBoundLamination}.

In Section \ref{NonExtensibilityQuestions}, we consider properties
of non-extensible combinatorial torus cube packings.
Firstly, we prove in Theorem \ref{ExtendibilityNM3} that combinatorial
cube packings with at least $2^n-3$ cubes are extensible to tilings.
In Propositions \ref{NumberParameters} and \ref{MinimalNrParam}, 
we prove that non-extensible combinatorial torus cube packings obtained
with strictly positive probability have at least $\frac{n(n+1)}{2}$
parameters
and that this number is attained by a combinatorial cube packing with $n+1$
cubes if $n$ is odd.
We conjecture that the number of parameters is at most $2^n-1$
(see Conjecture \ref{Power2conjecture}).
In Proposition \ref{6dimCases} we prove that in dimension $6$ the minimal
number of cubes in non-extensible combinatorial cube packings is $8$
and that none of those cube packings is attained with strictly positive
probability.
In Proposition \ref{StrangeDiscreteCases} we show that in dimension
$3$, $5$, $7$ and $9$, there exist combinatorial cube tilings obtained
with strictly positive probability and $\frac{n(n+1)}{2}$ parameters.

We now explain the origin of the model considered here.
Pal\'asti \cite{palasti} considered maximal packings
obtained from random packings
of cubes $[0,1]^n$ into $[0,x]^n$.
She conjectured that the expectation $E(M_x(n))$ of the packing
density $M_x(n)$ satisfies the limit 
\begin{equation}\label{FundamentalLimit}
\lim_{x\to\infty} \frac{E(M_x(n))}{x^n}=\beta_n\;.
\end{equation}
with $\beta_n=\beta_1^n$.
The value of $\beta_1$ is known since the work of R\'enyi \cite{renyi}
and in Penrose \cite{penrose01} the limit (\ref{FundamentalLimit}) is proved to
exist.
Note that based on simulations it is expected that $\beta_n > \beta_1^n$
and an experimental formula from simulations
\begin{equation*}
\beta_n^{1/n}-\beta_1 \simeq (n-1)(\beta_2^{1/2}-\beta_1)
\end{equation*}
is known \cite{solomon}.

The Itoh Ueda model \cite{ueda} is a variant of the above: one considers
packing of cubes $z+[0,2]^n$ with $z\in \ZZ^n$ into $[0,4]^n$.
It is proved in \cite{dip,poyarkov2,poyarkov} that the average number of
cubes satisfies the inequality
$E(M_2^C(n))\geq (\frac{3}{2})^n$ and some computer estimate of the average
density $\frac{1}{2^n} E(M_2^C(n))$ were obtained in \cite{yoshiaki}.
In \cite{cubetiling}, we considered the torus case, similar questions to
the one of this paper and a measure of regularity called second moment,
which has no equivalent here.

\section{Combinatorial cube packings}\label{SecCombCubPacking}

If $z+[0,1]^n\subset [0,2]^n$ and $z=(z_1, \dots,z_n)\in \frac{1}{N} \ZZ^n$
then $z_i\in \{0, \frac{1}{N}, \dots, 1\}$.
Take a discrete cube packing ${\mathcal{CP}}=(z^i+[0,1]^n)_{1\leq i\leq m}$ of $[0,2]^n$.
For a given coordinate $1\leq j\leq n$ we set $\phi(z^i_j)=t_{i,j}$ with
$t_{i,j}$ a parameter if $0<z^i_j<1$ and $\phi(z^i_j)=z^{i}_j$ if $z^{i}_j=0$
or $1$.
If $z^i=(z^i_1,z^i_2,\dots, z^i_n)$ then we set $\phi(z^i)=(\phi(z^i_1), \dots, 
\phi(z^i_n))$ and to ${\mathcal{CP}}$ we associate the combinatorial cube packing $\phi(\mathcal{CP})=(\phi(z^i)+[0,1]^n)_{1\leq i\leq m}$.

Take a torus discrete cube packing ${\mathcal{CP}}=(z^i+[0,1]^n)_{1\leq i\leq m}$ with $z^i\in \frac{1}{N}\ZZ^n$.
For a given coordinate $1\leq j\leq n$ we set $\phi(z^i_j)=t_{k, j}$
if $z^i_j\equiv \frac{k}{N}\pmod 2$ and $\phi(z^i_j)=t_{k, j}+1$
if $z^i_j\equiv \frac{k}{N}+1\pmod 2$ with $t_{k,j}$ a parameter.
Similarly, we set $\phi(z^i)=(\phi(z^i_1), \dots, \phi(z^i_n))$ and 
we define $\phi(\mathcal{CP})=(\phi(z^i)+[0,1]^n)_{1\leq i\leq m}$
the associated torus combinatorial cube packing.

In the remainder of this paper we do not use the above parameters
but instead renumber them into $t_1$, \dots, $t_N$.
Without loss of generality, we will always assume that different
coordinates have different parameters.
Of course we can define combinatorial cube packing without using
to discrete cube packings.
In the cube case, the relevant cubes are of the form $z+[0,1]^n$
with $z_i=0$, $1$ or some parameter $t$.
In the torus case, the relevant cubes are of the form $z+[0,1]^n$
with $z_i=t$ or $t+1$ and $t$ a parameter.
Two cubes $z^i+[0,1]^n$
and $z^{i'}+[0,1]^n$ are non-overlapping if there exist a coordinate $j$
such that $z^i_j\equiv z^{i'}_j+1\pmod 2$. In the cube case this means
that $z^i_j=0$ or $1$ and $z^{i'}_j=1-z^i_j$. In the torus case this
means that $z^i_j$ depends on the same parameter, say $t$, $z^i_j$, $z^{i'}_j=t$
or $t+1$ and $z^i_j\not= z^{i'}_j$.
A combinatorial cube packing is then a family of such cubes with
any two of them being non-overlapping. Notions of tilings and extensibility
are defined as well. Moreover, a discrete cube packing is extensible
if and only if its associated combinatorial cube packing is extensible.
Denote by $m({\mathcal{CP}})$ the number of cubes of a combinatorial
cube packing ${\mathcal{CP}}$ and by $N({\mathcal{CP}})$ its number of
parameters.
Denote by $Comb^C(n)$, $Comb^T(n)$, the set of
combinatorial cube packings of $[0,2]^n$,
respectively $\QuotS{\RR^n}{2\ZZ^n}$.

Given two combinatorial cube packings ${\mathcal{CP}}$ and ${\mathcal{CP}}'$
(either on cube or torus), we say that ${\mathcal{CP}}'$
is a {\em subtype} of ${\mathcal{CP}}$ if after assigning the parameter
of ${\mathcal{CP}}$ to $0$, $1$, or some parameter of 
${\mathcal{CP}}'$, we get ${\mathcal{CP}}'$.
So, necessarily $m({\mathcal{CP}}') = m({\mathcal{CP}})$ and 
$N({\mathcal{CP}}')\leq N({\mathcal{CP}})$ but the
reverse implication is not true in general.
A combinatorial cube packing is said to be {\em maximal}
if it is not the subtype of any other combinatorial cube packing.
Necessarily, a combinatorial cube packing ${\mathcal{CP}}$ is a subtype
of at least one maximal combinatorial cube packing ${\mathcal{CP}}'$.

Given a combinatorial cube packing ${\mathcal{CP}}$
the number of discrete cube packings ${\mathcal{CP}}'$ such that
$\phi({\mathcal{CP}}')={\mathcal{CP}}$ is denoted by $Nb({\mathcal{CP}}, N)$.
In the cube case we have $Nb({\mathcal{CP}}, N)=(N-1)^{N({\mathcal{CP}})}$.
The torus case is more complex but it is still possible to write 
explicit formulas: denote by $N_j({\mathcal{CP}})$ the number 
of parameters which occurs in the $j$-th coordinate of ${\mathcal{CP}}$.
We then get:
\begin{equation}\label{TorusCaseFormula}
Nb({\mathcal{CP}}, N)=\Pi_{j=1}^{n} \Pi_{k=1}^{N_{j}({\mathcal{CP}})} (2N-2(k-1)).
\end{equation}
The asymptotic order of $Nb({\mathcal{CP}}, N)$ is $(2N)^{N({\mathcal{CP}})}$, 
which shows that $Nb({\mathcal{CP}}, N)>0$ for $N$ large enough.
More specifically, $N_j({\mathcal{CP}})\leq 2^n$ so $Nb({\mathcal{CP}}, N)>0$
if $N\geq 2^n$.
Note that it is possible to have ${\mathcal{CP}}'$ a subtype of
${\mathcal{CP}}$ and $Nb({\mathcal{CP}}', N) > Nb({\mathcal{CP}}, N)$
for small enough $N$.

Denote by $f^{T}_N(n)$ the minimal number of cubes of 
non-extensible discrete torus cube packings
$(z^{i}+[0,1]^n)_{1\leq i\leq m}$ with $z^i\in \frac{1}{N}\ZZ^n$.
Denote by $f^T_{\infty}(n)$ the minimal number of cubes of non-extensible
combinatorial torus cube packings.

\begin{proposition}
For $n\geq 1$ we have 
$\lim_{N\to\infty} f^T_N(n)=f^T_{\infty}(n)$.
\end{proposition}
\proof A discrete cube packing ${\mathcal{CP}}$ is extensible if and only
if $\phi({\mathcal{CP}})$ is extensible. Thus $f^T_{\infty}(n)\geq f^T_N(n)$.
Take ${\mathcal{CP}}$ a non-extensible combinatorial torus cube packing
with the minimal number of cubes. By Formula (\ref{TorusCaseFormula})
there exist $N_0$ such that for $N> N_0$ we have $Nb({\mathcal{CP}}, N)>0$.
The discrete cube packings ${\mathcal{CP}}'$ with $\phi({\mathcal{CP}}')={\mathcal{CP}}$ are non-extensible.
So, we have $\lim_{N\to\infty} f^T_N(n)=f^T_{\infty}(n)$. \qed

In the cube case we have for $N\geq 2$ the equality $f^{C}_{N}(n)=1$.

Two combinatorial cube packings ${\mathcal{CP}}$ and ${\mathcal{CP}}'$ are said
to be equivalent if after a renumbering of the coordinates, parameters and
cubes of ${\mathcal{CP}}$ one gets ${\mathcal{CP}}'$.
The automorphism group of a combinatorial cube packing is the group of
equivalences of ${\mathcal{CP}}$ preserving it.
Testing equivalences and computing stabilizers can be done using the
program {\tt nauty} \cite{nauty},
which is a graph theory program for testing whether two
graphs are isomorphic or not and computing the automorphism group.
The method is to associate to a given combinatorial cube packing
${\mathcal{CP}}$ a graph $Gr({\mathcal{CP}})$, which characterize
isomorphism and automorphisms.
The method used to find such a graph $Gr({\mathcal{CP}})$
are explained in the user manual of {\tt nauty}
and the corresponding programs are available from \cite{dutourPackProb}.

We now explain the sequential random cube packing.
Given a discrete cube packing ${\mathcal{CP}}=(z^{i}+[0,1]^n)_{1\leq i\leq m}$
denote by $Poss({\mathcal{CP}})$ the set of cubes $z+[0,1]^n$ with
$z\in \frac{1}{N}\ZZ^n$ which do not overlap with ${\mathcal{CP}}$.
Every possible cube $z+[0,1]^n$ is selected with equal probability
$\frac{1}{|Poss({\mathcal{CP}})|}$.
The sequential random cube packing process is thus a process that
add cubes until 
the discrete cube packing is non-extensible or is a tiling.

Fix a combinatorial cube packing ${\mathcal{CP}}$, $N\geq 2^n$
and a discrete cube packing ${\mathcal{CP}}'$ such that
$\phi({\mathcal{CP}}')={\mathcal{CP}}$. 
To any cube $w+[0,1]^n\in Poss({\mathcal{CP}}')$
we associate the combinatorial cube packing ${\mathcal{CP}}_w=\phi({\mathcal{CP}}'\cup \{w+[0,1]^n\})$.
The set $Poss({\mathcal{CP}}')$ is partitioned 
into classes $Cl_1$, \dots, $Cl_r$ with two cubes 
$w+[0,1]^n$ and $w'+[0,1]^n$ in the same class if 
${\mathcal{CP}}_{w}={\mathcal{CP}}_{w'}$.
The combinatorial cube packing associated to $Cl_i$ is
denoted by ${\mathcal{CP}}_i$.
The set $\{{\mathcal{CP}}_1, \dots, {\mathcal{CP}}_r\}$ of classes
depends only on ${\mathcal{CP}}$. If we had chosen some $N\leq 2^n$, then
some of the preceding ${\mathcal{CP}}_i$ might not have occurred.
So, we have
\begin{equation*}
|Cl_i(N)|=\frac{Nb({\mathcal{CP}}_i, N)}{Nb({\mathcal{CP}}, N)}
\end{equation*}
and we can define the probability $p({\mathcal{CP}}, {\mathcal{CP}}_i, N)$ of obtaining a discrete cube packing of combinatorial type ${\mathcal{CP}}_i$ from a discrete cube packing of combinatorial type ${\mathcal{CP}}$:
\begin{equation*}
p({\mathcal{CP}}, {\mathcal{CP}}_i, N)
=\frac{|Cl_i(N)|}{|Cl_1(N)|+\dots+|Cl_r(N)|}
=\frac{Nb({\mathcal{CP}}_i, N)}{Nb({\mathcal{CP}}_1, N)+\dots+Nb({\mathcal{CP}}_r, N)} .
\end{equation*}
Given a combinatorial cube packing ${\mathcal{CP}}$ with $m$ cubes a 
{\em path} $p=\{{\mathcal{CP}}^0, {\mathcal{CP}}^1, \dots, {\mathcal{CP}}^m\}$
is a way of obtaining ${\mathcal{CP}}$ by adding one cube at a time
starting from
${\mathcal{CP}}^0=\emptyset$ and ending at ${\mathcal{CP}}^m={\mathcal{CP}}$.
The probability to obtain ${\mathcal{CP}}$ along a path $p$ is
\begin{equation*}
p({\mathcal{CP}}, p, N)=
p({\mathcal{CP}}^0, {\mathcal{CP}}^1, N)
\times\dots\times
p({\mathcal{CP}}^{m-1}, {\mathcal{CP}}^m, N).
\end{equation*}
The probability $p({\mathcal{CP}}, N)$ to obtain ${\mathcal{CP}}$
is the sum over all the paths $p$ leading to ${\mathcal{CP}}$
of $p({\mathcal{CP}}, p, N)$. The probabilities $p({\mathcal{CP}}, p, N)$
and $p({\mathcal{CP}}, N)$ are fractional functions of $N$, which
implies that the limit value $p({\mathcal{CP}}, \infty)$, 
$p({\mathcal{CP}}, p, \infty)$ and $p({\mathcal{CP}}, {\mathcal{CP}}', \infty)$
are well defined.

As $N$ goes to $\infty$ we have the asymptotic behavior
\begin{equation*}
|Cl_i(N)|\simeq (2N)^{nb_i}
\end{equation*}
with $nb_i=N({\mathcal{CP}}_i)-N({\mathcal{CP}})$ the number
of new parameters in ${\mathcal{CP}}_i$ as compared
with ${\mathcal{CP}}$.
Clearly as $N$ goes to $\infty$ only the classes with the largest $nb_i$
have $p({\mathcal{CP}}, {\mathcal{CP}}_i, \infty)>0$.
If $Cl_i$ is such a class then we get 
\begin{equation*}
p({\mathcal{CP}}, {\mathcal{CP}}_i, \infty)=\frac{1}{r'}
\end{equation*}
with $r'$ the number of classes $Cl_i$ having the largest $nb_i$ and otherwise 
$p({\mathcal{CP}}, {\mathcal{CP}}_i, \infty)=0$.
Analogously, for a path $p$ leading to ${\mathcal{CP}}$ we can define
$p({\mathcal{CP}}, p, \infty)$ and $p({\mathcal{CP}}, \infty)$.
We say that a combinatorial cube packing ${\mathcal{CP}}$ is obtained
with {\em strictly positive probability} if $p({\mathcal{CP}}, \infty)>0$ 
that is for at least one path $p$ we have $p({\mathcal{CP}}, p, \infty)>0$.
For a path $p=\{{\mathcal{CP}}^0, {\mathcal{CP}}^1, \dots, {\mathcal{CP}}^m\}$
we have $p({\mathcal{CP}}, p, \infty)>0$ if and only if every
${\mathcal{CP}}^i$ has $N({\mathcal{CP}}^i)$ maximal among all possible
extensions from ${\mathcal{CP}}^{i-1}$. This implies that each 
${\mathcal{CP}}^i$ is maximal, i.e. is not the subtype of another type.
As a consequence, we can define a sequential random cube packing process for
combinatorial cube packing ${\mathcal{CP}}$ obtained with strictly positive
probability and compute their probability $p({\mathcal{CP}}, \infty)$.

A combinatorial cube packing ${\mathcal{CP}}$ is said to have order
$k=\tord({\mathcal{CP}})$ if $p({\mathcal{CP}}, N)=\frac{1}{(N-1)^k}f(N)$
with $\lim_{N\to\infty} f(N)\in \RR_+^*$.
A combinatorial cube packing is of order $0$ if and only if it is obtained
with strictly positive probability.

Let us denote by $M_N^C(n)$, $M_N^T(n)$ the random variables of number
of cubes of those non-extensible cube packings and by $E(M_N^C(n))$,
$E(M_N^T(n))$ their expectation.
From the preceding
discussion we have
\begin{equation*}
E(M_N^U(n))=\sum_{ {\mathcal{CP}}\in Comb^U(n)} p({\mathcal{CP}}, N) m({\mathcal{CP}})
\mbox{~~with~~}U\in \{C, T\}.
\end{equation*}
Denote by $f^T_{>0,\infty}(n)$ the minimal number of cubes of 
non-extensible combinatorial torus cube packings obtained with
strictly positive probability.

\begin{figure}
\begin{center}
\begin{minipage}[b]{7cm}
\centering
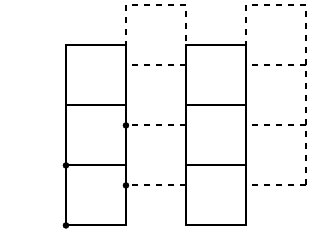\par
A combinatorial cube tiling obtained with probability $\frac{1}{2}$.
\end{minipage}
\begin{minipage}[b]{7cm}
\centering
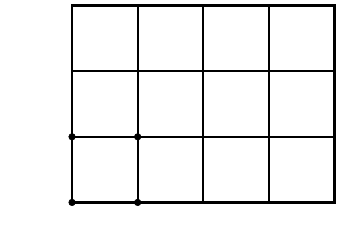\par
A combinatorial cube tiling obtained with probability $0$.
\end{minipage}
\end{center}
\caption{Two $2$-dimensional torus combinatorial cube tilings}
\label{2dimContinuousCubePacking}
\end{figure}

\begin{figure}
\begin{center}
\begin{minipage}[b]{3.3cm}
\centering
\epsfig{file=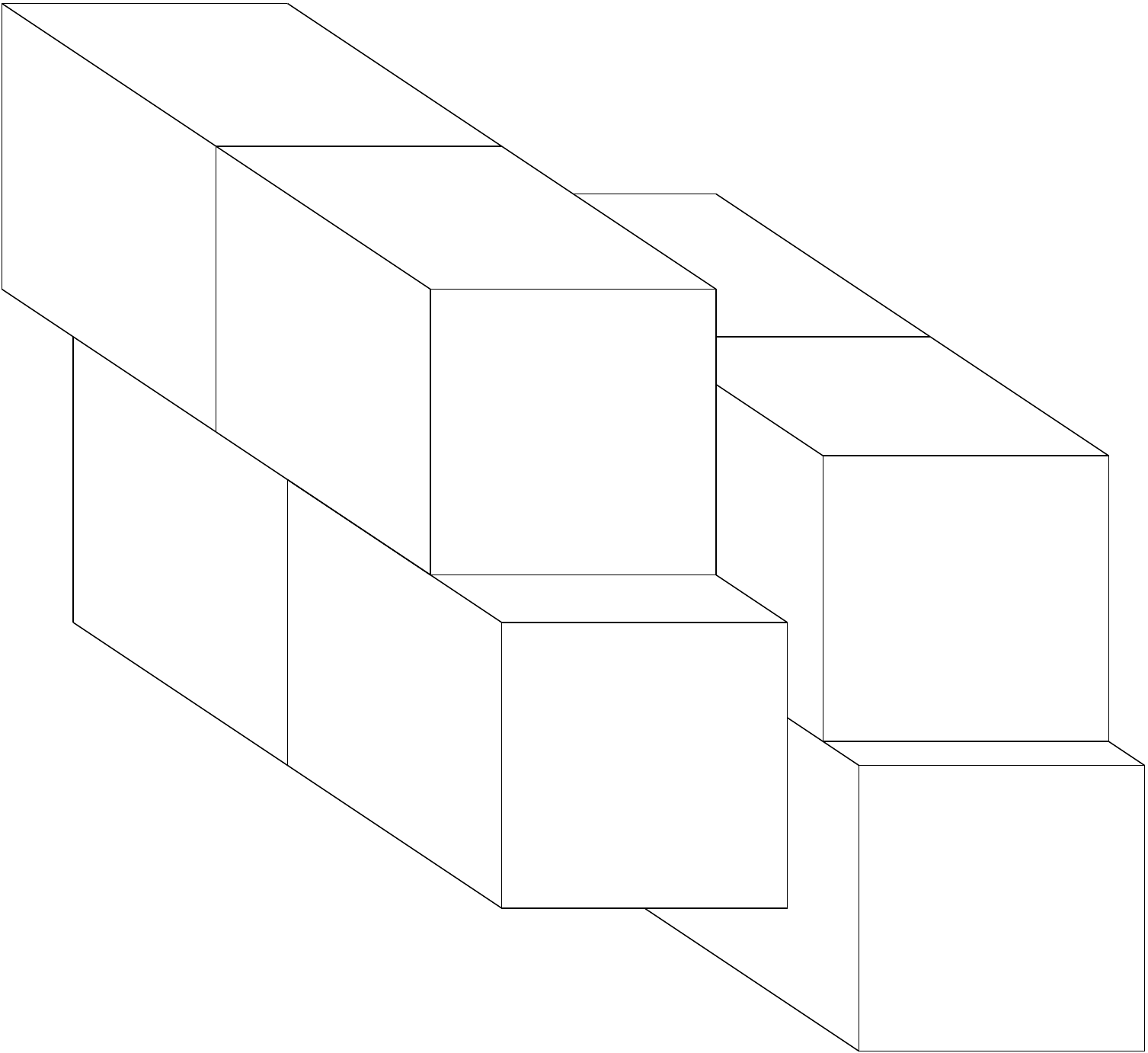, width=28mm}\par
$7$ parameters,\par
probability $\frac{1}{3}$
\end{minipage}
\begin{minipage}[b]{3.3cm}
\centering
\epsfig{file=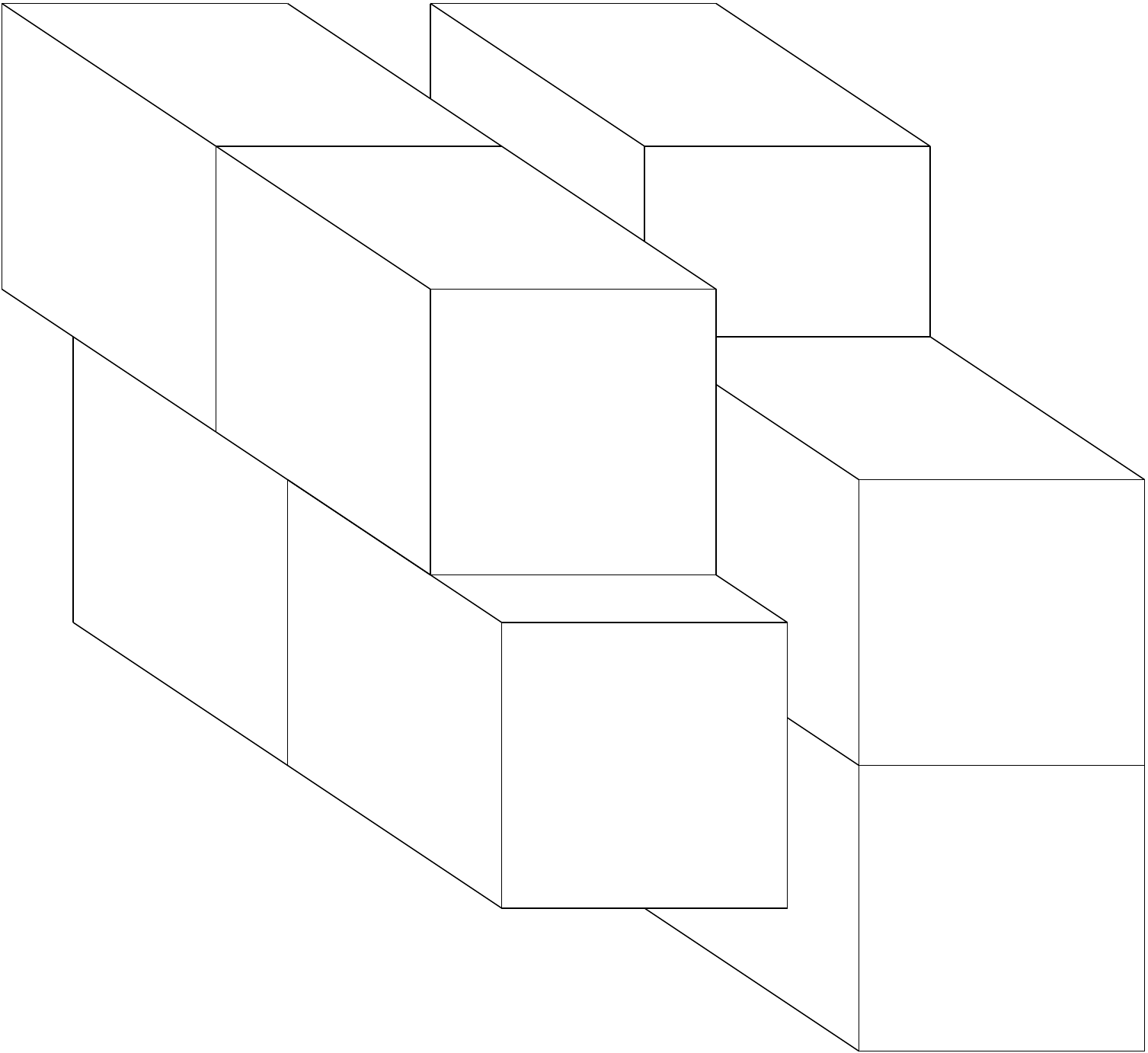, width=28mm}\par
$7$ parameters,\par
probability $\frac{1}{3}$
\end{minipage}
\begin{minipage}[b]{3.3cm}
\centering
\epsfig{file=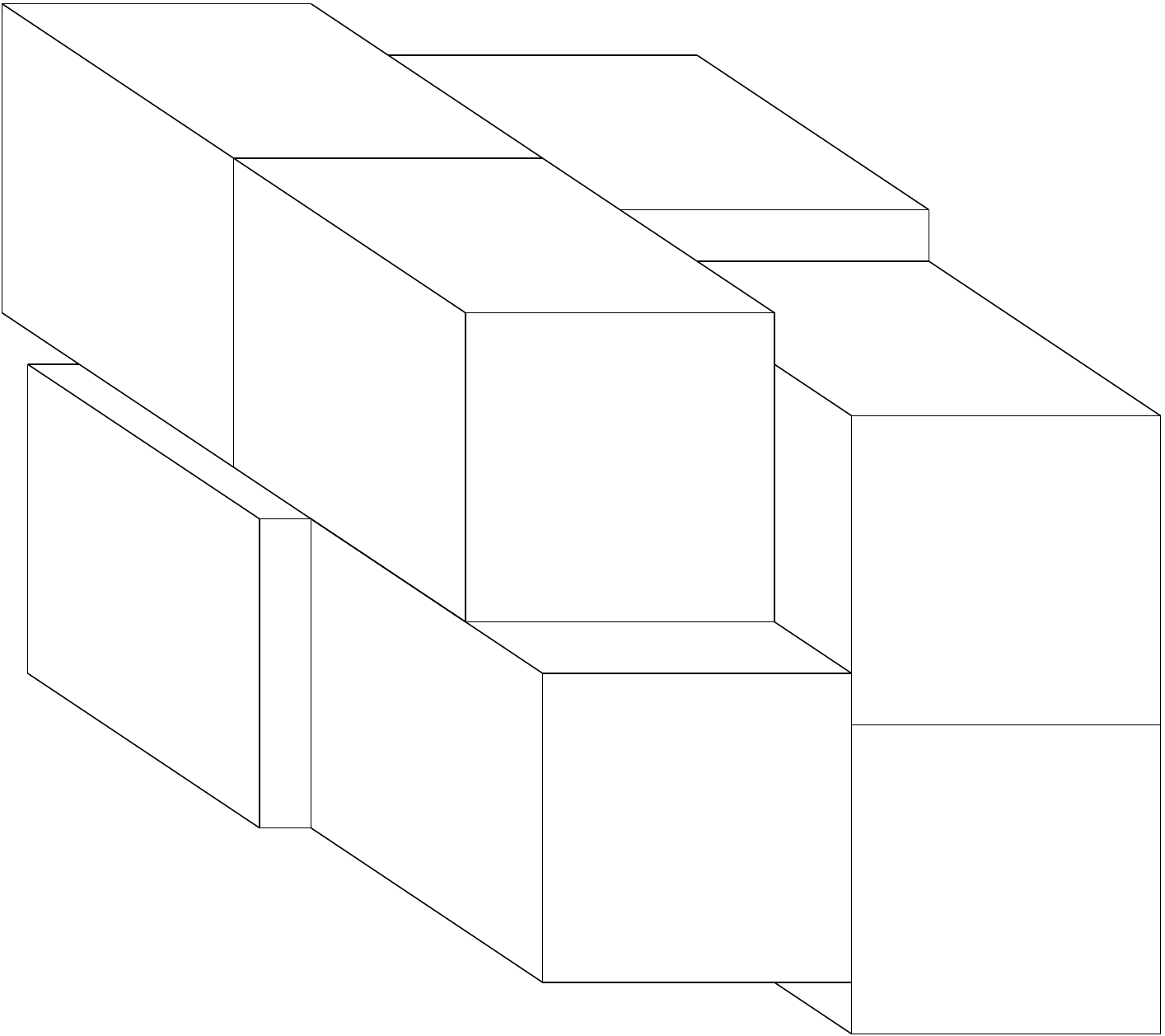, width=28mm}\par
$6$ parameters,
probability $\frac{5}{18}$
\end{minipage}
\begin{minipage}[b]{3.3cm}
\centering
\includegraphics[width=32mm,bb=109 241 457 591,clip]{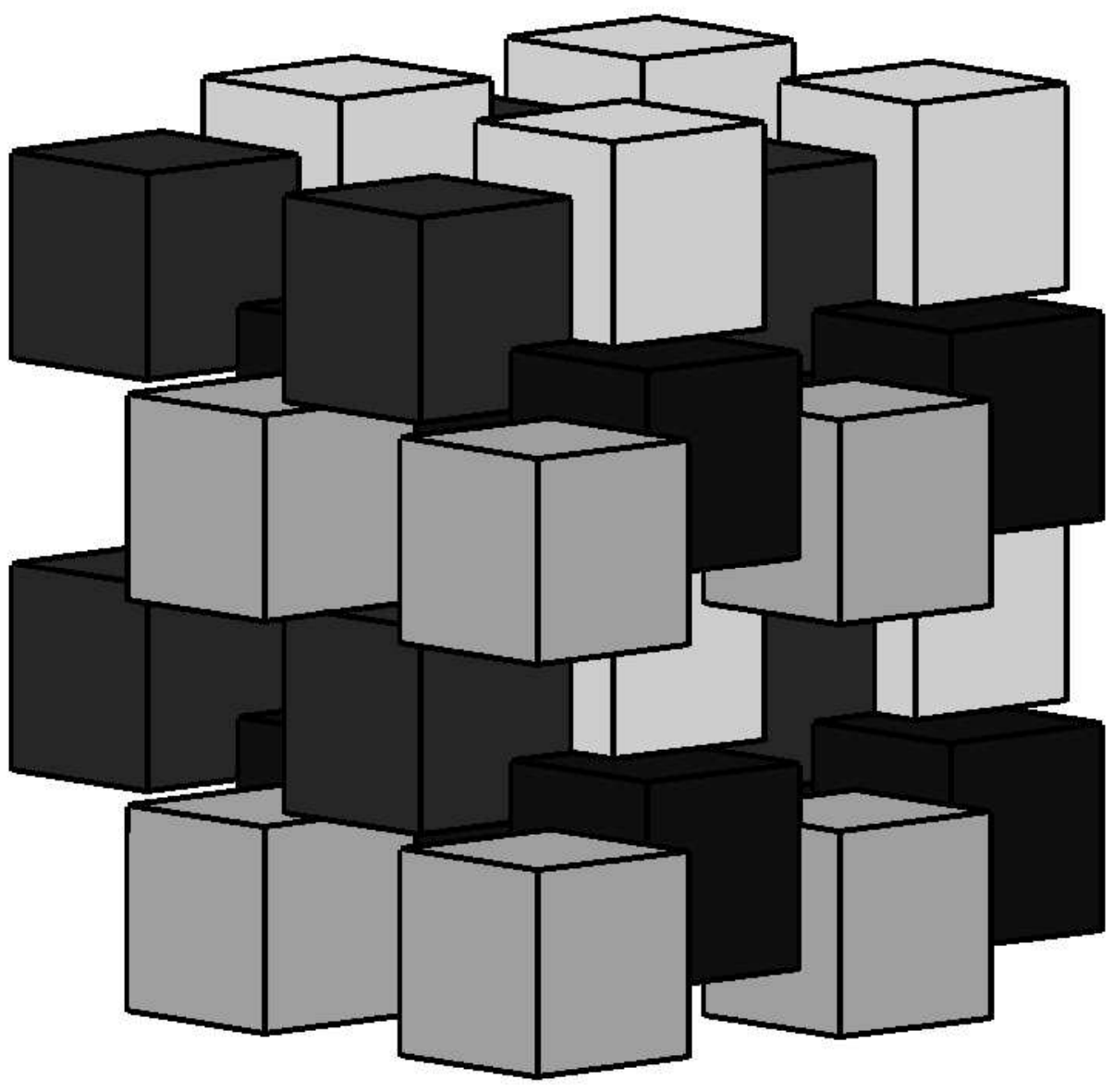}\par
$6$ parameters,
probability $\frac{1}{18}$
\end{minipage}

\end{center}
\caption{The $3$-dimensional combinatorial cube packings obtained with strictly positive probability; two laminations over $2$-dimensional cube tilings, the rod tiling and the smallest non-extensible cube packing}
\label{EnumerationDim3}
\end{figure}

In dimension $2$ (see Figure \ref{2dimContinuousCubePacking}), there
are three combinatorial cube tilings. 
One of them is attained with probability $0$; it is a subtype of
the remaining two which
are equivalent and attained with probability $\frac{1}{2}$.
By applying the random cube packing process and doing reduction
by isomorphism, one obtains the $3$-dimensional combinatorial cube packings
obtained with strictly positive probability (see Figure \ref{EnumerationDim3}).
The non-extensible cube packing shown on this figure already occurs
in \cite{lagariasEmail,cubetiling}.
In dimension $4$, the same enumeration method works
(see Table \ref{DataInformation})
but dimension $5$ is computationally too difficult to enumerate.

\begin{table}
\begin{center}
\begin{tabular}{|c|c|c|c|c|c|c|}
\hline
      &$n$                   & $1$  & $2$ & $3$ & $4$ & $5$\\
\hline
\hline
$N=\infty$ &Nr cube tilings  & $1$  & $1$ & $3$ & $32$& ?\\
\hline
      &Nr non-extensible       & $0$  & $0$ & $1$ & $31$& ?\\
      &cube packings           &      &     &     &     &\\
\hline
      & $f^T_{>0,\infty}(n)$ & $2$  & $4$ & $4$ & $6$ & $6$\\\hline
      & $\frac{1}{2^n}E(M_{\infty}^T(n))$   & $1$  & $1$ & $\frac{35}{36}$ & $\frac{15258791833}{16102195200}$ & ?\\
\hline
\hline
$N=2$ &Nr cube tilings  & 1  & 2 & 8 & 744 & ?\\
\hline
      &Nr cube packings & 0  & 0 & 1 & 139 & ?\\
\hline
      &$f^T_2(n)$         & 2  & 4 & 4 & 8 & $10\leq f^T_2(5)\leq 12$\\
\hline
\end{tabular}
\end{center}
\caption{Number of packings and tilings for the case $N=\infty$
and $N=2$ (see \cite{cubetiling})}
\label{DataInformation}
\end{table}

\section{Discrete Random cube packings of the cube}\label{DiscreteCubePackCube}
We compute here the polynomials $C_k(n)$ occurring in Equation (\ref{ExpansionCombCubPack})
for $k\leq 6$. We compute the first three polynomials by an elementary method.

\begin{lemma}\label{ToyLemma}
Put the cube $z^1+[0,1]^n$ in $[0,2]^n$ and write
\begin{equation*}
I=\{i\;\;:\;\;z^1_i=0\mbox{~or~}N\}
\end{equation*}
then do sequential random discrete cube packing.
\begin{itemize}
\item[(i)] The minimal number of cubes in the packing is $|I|+1$.
\item[(ii)] The expected number of cubes in the packing
is $|I|+1+O\left(\frac{1}{N+1}\right)$.
\end{itemize}

\end{lemma}
\proof Let us prove (i).
If $|I|=0$, then clearly one cannot insert any more cubes.
We will assume $|I|>1$ and do a reasoning by induction on $|I|$.
If one puts another cube $z^2+[0,1]^n$, there should exist
an index $i\in I$ such that
$|z^2_i-z^1_i|=1$. Take an index $j\not= i$ such that $z^2_j$ is $0$ or $1$.
The set of possibilities to add a subsequent cube is larger if
$z_j\in \{0,1\}$ than if $z_j\in \{\frac{1}{N},\dots,\frac{N-1}{N} \}$.
So, one can assume that for $j\not= i$, one has $0 < z^2_j < 1$.
This means that any cube $z+[0,1]^n$ in subsequent insertion should
satisfy $|z_i-z^2_i|=1$, i.e. $z_i=z^1_i$. So, the sequential
random cube packing can be done in one dimension less, starting with
${z'}^1=(z_1,\dots, z_{i-1},z_{i+1},\dots, z_n)$.
The induction hypothesis applies.
Assertion (ii) follows easily by looking at the above
process. For a given $i$ the choice of $z^2$ with $0< z^2_j < 1$ for
$j\not= i$ is the one with probability $1-O\left(\frac{1}{N+1}\right)$.
So, all neglected possibilities have probability $O\left(\frac{1}{N+1}\right)$
and with probability $1-O\left(\frac{1}{N+1}\right)$ the number of cubes is the
minimal possible. \qed

See below the $2$-dimensional possibilities:

\begin{center}
\begin{minipage}{3.3cm}
\centering
\resizebox{1.5cm}{!}{\includegraphics{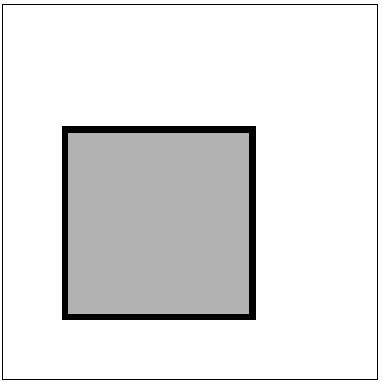}}\par
$|I|=0$
\end{minipage}
\begin{minipage}{3.3cm}
\centering
\resizebox{1.5cm}{!}{\includegraphics{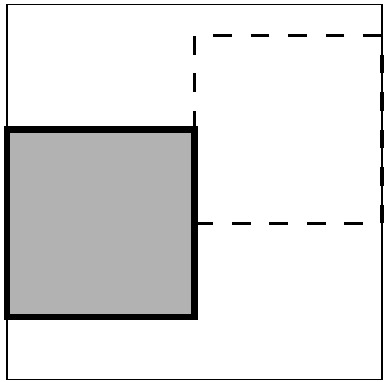}}\par
$|I|=1$
\end{minipage}
\begin{minipage}{3.3cm}
\centering
\resizebox{1.5cm}{!}{\includegraphics{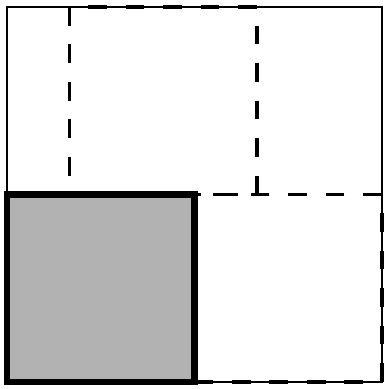}}\par
$|I|=2$
\end{minipage}
\end{center}

The random variable $M_N^C(n)$ is the number of cubes in the obtained
non-extensible cube-packing.
$E(M_N^C(n))$ is the expected number of cubes and $E(M_N^C(n)~|~k)$ the expected
number of cubes obtained by imposing the condition that the first
cube $z^1+[0,1]^n$ has $|\{i\;\;:\;\;z^1_i=0\mbox{~or~}1\}|=k$.

\begin{theorem}\label{SecondExponent}
For any $n\geq 1$, we have
\begin{equation*}
E(M_N^C(n))=1+\frac{2n}{N+1}+\frac{4n(n-1)}{(N+1)^2}+O\left(\frac{1}{N+1}\right)^3\mbox{~as~}N\to\infty.
\end{equation*}
\end{theorem}
\proof If one chooses a vector $z$ in $\{0,\dots,N\}^{n}$ the probability
that $|\{i\;\;:\;\;z^1_i=0\mbox{~or~}N\}|=k$ is $\left\lbrace\frac{2}{N+1}\right\rbrace^{k} \left\lbrace\frac{N-1}{N+1}\right\rbrace^{n-k} {n \choose k}$.
Conditioning over $k\in \{0, 1, \dots, n\}$, one obtains
\begin{equation}\label{ConditioningFormula}
E(M_N^C(n))=\sum_{k=0}^n \left\lbrace\frac{2}{N+1}\right\rbrace^{k} \left\lbrace\frac{N-1}{N+1}\right\rbrace^{n-k} {n \choose k} E(M_N^C(n)~|~k).
\end{equation}
So one gets
\begin{equation*}
\begin{array}{rcl}
E(M_N^C(n))&=&\left(\frac{N-1}{N+1}\right)^n  E(M_N^C(n)~|~0)+n\frac{2}{N+1} \left(\frac{N-1}{N+1}\right)^{n-1} E(M_N^C(n)~|~1)\\
&+&n(n-1)\frac{2}{(N+1)^2} \left(\frac{N-1}{N+1}\right)^{n-2} E(M_N^C(n)~|~2)+O\left(\frac{1}{N+1}\right)^3.
\end{array}
\end{equation*}
Clearly $E(M_N^C(n)~|~0)=1$ and $E(M_N^C(n)~|~1)=1+E(M_N^C(n-1))$.
By Lemma \ref{ToyLemma}, $E(M_N^C(n)~|~2)=3+O\left(\frac{1}{N+1}\right)$.
Then, one has $E(M_N^C(n))=1+O\left(\frac{1}{N+1}\right)$ and 
\begin{equation*}
\begin{array}{rcl}
E(M_N^C(n))&=&\left(1-\frac{2}{N+1}\right)^n+\frac{2n}{N+1}\left(1-\frac{2}{N+1}\right)^{n-1}(2+O\left(\frac{1}{N+1}\right))\\
&+&O\left(\frac{1}{(N+1)^2}\right)\\
&=&\left\{1-\frac{2n}{N+1}+O\left(\frac{1}{(N+1)^2}\right)\right\}+\frac{4n}{N+1}+O\left(\frac{1}{(N+1)^2}\right)\\
&=&1+\frac{2n}{N+1}+O\left(\frac{1}{(N+1)^2}\right).
\end{array}
\end{equation*}
Inserting this expression into $E(M_N^C(n))$ and Formula (\ref{ConditioningFormula}) one gets the result. \qed

So, we get $C_0(n)=1$, $C_1(n)=2n$ and $C_2(n)=4n(n-2)$.
In order to compute $C_k(n)$ in general we use methods similar to the
ones of Section \ref{SecCombCubPacking}.
Given a cube $z+[0,1]^n$ with $z_i\in \{0, \frac{1}{N}, \dots, 1\}$ we define
a face of the cube $[0,1]^n$ in the following way: if $z_i=0$ or $1$ then we set
$\psi(z_i)=0$ or $1$ whereas if $0<z_i<1$ we set $\psi(z_i)=t_i$ with $t_i$
a parameter. 
When the parameters $t_i$ of the vector $(\psi(z_1), \dots, \psi(z_n))$
vary in $]0,1[$ this vector describes a face of the cube $[0,1]^n$, which
we denote by $\psi(z)$.
This construction was presented for the first time
in \cite{poyarkov,poyarkov2}.

If $F$ and $F'$ are two faces of $[0,1]^n$, then we say that
$F$ is a sub-face of $F'$ and write $F\subset F'$ if $F$ is
included in the closure of $F'$.
A subcomplex of the hypercube $[0,1]^n$ is a set of faces,
which contains all its sub-faces.
If ${\mathcal{CP}}$ is a cube packing in $[0,2]^n$, then the
vectors $z$ such that $z+[0,1]^n$ is a cube which we can add
to it are indexed by the faces
of a subcomplex $[0,1]^n$ with the dimension giving the exponent of
$(N-1)^k$.
The dimension of a complex is the highest dimension of its faces.
Given a discrete cube packing ${\mathcal{CP}}$, we have seen in
Section \ref{SecCombCubPacking} that the size of $Poss({\mathcal{CP}})$
depends only on the combinatorial type $\phi({\mathcal{CP}})$.
In the cube case which we
consider in this section $Poss({\mathcal{CP}})$ itself depends only
on the combinatorial type.

\begin{theorem}\label{MainTheoremRigidBoundaryCase}
There exist polynomials $C_{k}(n)$ of $n$ with $\tdeg\,C_{k}=k$
such that for any $n$ and $N>1$ one has:
\begin{equation*}
E(M_N^C(n))= \sum_{k=0}^{\infty} \frac{C_{k}(n)}{(N-1)^k}.
\end{equation*}
The polynomials $C_k(n)$ are given in Table \ref{PolynomesRkn}.
\end{theorem}
\proof The image $\psi(Poss({\mathcal{CP}}))$ is an
union of faces of $[0,1]^n$, i.e. a subcomplex of the complex $[0,1]^n$.
Denote by $\tdim(F)$ the dimension of a
face $F$ of the cube $[0,1]^n$.
Denote by $Poss(F)$ the set of vectors $z\in \{0, \frac{1}{N},\dots, 1\}^n$
with $\psi(z)=F$.
we have the formula:
\begin{equation*}
|Poss(F)|=(N-1)^{\tdim(F)}\mbox{~~and~~}|Poss({\mathcal{CP}})|=\sum_{F} (N-1)^{\tdim(F)}.
\end{equation*}
The cubes, whose corresponding face in $[0,1]^n$ have dimension
$\tdim(\psi(Poss({\mathcal{CP}})))$ have the highest probability
of being obtained.
If one seeks the expansion of $E(M_N^C(n))$
up to order $k$ and if ${\mathcal{CP}}$ is 
of order $\tord({\mathcal{CP}})$ then we need to compute the
faces of $\psi(Poss({\mathcal{CP}}))$ of dimension at
least $\tdim(\psi(Poss({\mathcal{CP}}))) - (k-\tord({\mathcal{CP}}))$.
The probabilities are then obtained in the following way:
\begin{equation}\label{probaPFN}
p(F, N)=\frac{(N-1)^{\tdim(F)}}{\sum_{F'\in \psi(Poss({\mathcal{CP}}))\mbox{~with~}\tdim(F')\geq \dim(\psi(Poss({\mathcal{CP}}))) - (k-\tord({\mathcal{CP}}))} (N-1)^{\tdim(F')}}.
\end{equation}
The enumeration algorithm is then the following:
\begin{flushleft}
\smallskip
\textbf{Input:} Exponent $k$.\\
\textbf{Output:} List ${\mathcal L}$ of all inequivalent combinatorial types
of non-extensible cube packings ${\mathcal{CP}}$ with order
at most $k$ and their probabilities $p({\mathcal{CP}}, N)$ with an error of $O\left(\frac{1}{(N+1)^{k+1}}\right)$.\\
\smallskip
${\mathcal T} \leftarrow \{\emptyset\}$.\\
${\mathcal L} \leftarrow \emptyset$\\
\textbf{while} there is a ${\mathcal{CP}} \in {\mathcal T}$ \textbf{do}\\
\hspace{2ex} ${\mathcal T} \leftarrow {\mathcal T} \setminus \{{\mathcal{CP}}\}$\\
\hspace{2ex} $\psi(Poss({\mathcal{CP}})) \leftarrow$ the complex of all possibilities of adding a cube to ${\mathcal{CP}}$\\
\hspace{2ex} ${\mathcal F} \leftarrow$ the faces of $\psi(Poss({\mathcal{CP}}))$ of dimension at least\\
\hspace{8ex} $\tdim(\psi(Poss({\mathcal{CP}}))) - (k-\tord({\mathcal{CP}}))$\\
\hspace{2ex} \textbf{if} ${\mathcal F} = \emptyset$ \textbf{then}\\
\hspace{2ex} \hspace{2ex} \textbf{if} ${\mathcal{CP}}$ is equivalent to a ${\mathcal{CP}}'$ in ${\mathcal L}$ \textbf{then}\\
\hspace{2ex} \hspace{2ex} \hspace{2ex} $p({\mathcal{CP}}', N) \leftarrow p({\mathcal{CP}}', N)+p({\mathcal{CP}}, N)$\\
\hspace{2ex} \hspace{2ex} \textbf{else}\\
\hspace{2ex} \hspace{2ex} \hspace{2ex} ${\mathcal L}\leftarrow {\mathcal L}\cup \{{\mathcal{CP}}\}$\\
\hspace{2ex} \hspace{2ex} \textbf{end if}\\
\hspace{2ex} \textbf{else}\\
\hspace{2ex} \hspace{2ex} \textbf{for} $C \in {\mathcal F}$ \textbf{do}\\
\hspace{2ex} \hspace{2ex} \hspace{2ex} ${\mathcal{CP}}_{new} \leftarrow {\mathcal{CP}}\cup \{C\}$\\
\hspace{2ex} \hspace{2ex} \hspace{2ex} $p({\mathcal{CP}}_{new}, N) \leftarrow p({\mathcal{CP}}, N) p(C, N)$\\
\hspace{2ex} \hspace{2ex} \hspace{2ex} \textbf{if} ${\mathcal{CP}}_{new}$ is equivalent to a ${\mathcal{CP}}'$ in ${\mathcal T}$ \textbf{then}\\
\hspace{2ex} \hspace{2ex} \hspace{2ex} \hspace{2ex} $p({\mathcal{CP}}', N) \leftarrow p({\mathcal{CP}}', N)+p({\mathcal{CP}}_{new}, N)$\\
\hspace{2ex} \hspace{2ex} \hspace{2ex} \textbf{else}\\
\hspace{2ex} \hspace{2ex} \hspace{2ex} \hspace{2ex} ${\mathcal T}\leftarrow {\mathcal T} \cup \{{\mathcal{CP}}_{new}\}$\\
\hspace{2ex} \hspace{2ex} \hspace{2ex} \textbf{end if}\\
\hspace{2ex} \hspace{2ex} \textbf{end for}\\
\hspace{2ex} \textbf{end if}\\
\textbf{end while}\\
\end{flushleft}
Let us prove that the coefficients $C_k(n)$ are polynomials in the
dimension $n$.
If $C$ is the cube $[0,1]^n$ then the number of faces of codimension 
$l$ is $2^{l} {n \choose l}$, i.e. a polynomial in $n$ of degree $l$.
Suppose that a cube packing ${\mathcal{CP}}=(z^i+[0,1]^n)_{1\leq i\leq m}$
has $0< z^i_j < 1$ for $n'\leq j\leq n$.
Then all faces $F$ of $\psi(Poss({\mathcal{CP}}))$ of maximal dimension
$d=\dim(\psi(Poss({\mathcal{CP}})))$ have $0 < z_j < 1$ for $n'\leq j\leq n$
and $z\in F$.
When one chooses a subface of $F$ of dimension $d-l$, we have to choose
some coordinates $j$ to be equal to $0$ or $1$.
Denote by $l'$ the number of such coordinates $j$ with $n'\leq j\leq n$.
There are $2^{l'} {n+1-n' \choose l'}$ choices and they are all equivalent.
There are still $l-l'$ choices to be made for $j\leq n'-1$ but this number
is finite so in all cases the faces of $\psi(Poss({\mathcal{CP}}))$ 
of dimension at least $d-l$ can
be grouped in a finite number of classes with the size of the classes depending
on $n$ polynomially.
Moreover, the number of classes of dimension $d$ is finite so the 
term of higher order in the denominator of Equation (\ref{probaPFN})
is constant and the coefficients of the expansion of $p(F,N)$
are polynomial in $n$. \qed

\begin{table}
\begin{center}
\begin{tabular}{|c|c|c|}
\hline
$k$   &  $|Comb^C_k|$   & $C_{k}(n)$\\
\hline
$0$   &  $1$        & $1$\\
$1$   &  $2$        & $2n$\\
$2$   &  $3$        & $4n(n-2)$\\
$3$   &  $7$        & $\frac{1}{3}\{28n^3-153n^2+149n\}$\\
$4$   &  $18$       & $\frac{1}{2\cdot 3^2\cdot 5}\{2016n^4-21436n^3+58701n^2-40721n\}$\\
$5$   &  $86$       & $\frac{1}{2^2 3^3 \cdot 5\cdot 7}\{208724n^5-3516724n^4+18627854n^3-35643809n^2+20444915n\}$\\
$6$   &  $1980$     & $\frac{1}{2^8 3^{11} 5^5 7^3 11^3 \cdot 13\cdot 17}\{1929868729224214329703n^6-46928283796201160537385 n^5$\\
      &             &$+397379056595496330171955 n^4-1442659974291080413770375 n^3$\\
      &             &$+2205275555952621337847422 n^2-1115911322466787143241320 n\}$\\\hline
\end{tabular}
\end{center}
\caption{The polynomials $C_{k}(n)$. $|Comb^C_k|$ is the number of types of
combinatorial cube packings ${\mathcal{CP}}$ with $\tord({\mathcal{CP}})\leq k$}
\label{PolynomesRkn}
\end{table}

\section{Combinatorial torus cube packings and lamination construction}\label{LaminationTorusCubePacking}

\begin{lemma}\label{TramLemma}
Let ${\mathcal{CP}}$ be a non-extensible combinatorial torus cube packing.

(i) Every parameter $t$ of ${\mathcal{CP}}$ occurs, which occurs as $t$ also occurs as $t+1$.

(ii) Let $C_1, \dots, C_k$ be cubes of ${\mathcal{CP}}$ and $C$
a cube which does not overlap with ${\mathcal{CP}}'={\mathcal{CP}}-\{C_1, \dots, C_k\}$.
The number of parameters of $C$, which does not occur in ${\mathcal{CP}}'$
is at most $k-1$.

\end{lemma}
\proof (i) Suppose that a parameter $t$ of ${\mathcal{CP}}$ occurs as $t$ but not as $t+1$ in the coordinates of the cubes.
Let $C=z+[0,1]^n$ be a cube having $t$ in its $j$-th coordinate.
If $C'=z'+[0,1]^n$ is a cube of ${\mathcal{CP}}$, then there exist a
coordinate $j'$ such that $z'_{j'}\equiv z_{j'}+1\pmod 2$.
Necessarily $j'\not= j$ since $t+1$ does not occur, so $C+e_j$ does not
overlap with $C'$ as well and obviously $C+e_j$ does not overlap with $C$.

(ii) Let $C=z+[0,1]^n$ be a cube which does not overlap with the cubes of 
${\mathcal{CP}}'$. 
Suppose that $z$ has $k$ coordinates $i_1< \dots < i_k$ such that their
parameters $t_1, \dots, t_k$ do not occur in ${\mathcal{CP}}'$.
If $C_j=z^j+[0,1]^n$, then we fix $z_{i_j}\equiv z^{j}_{i_j}+1\pmod 2$ for $1\leq j\leq k$ so that $C$ does not overlap with ${\mathcal{CP}}$.
This contradicts the fact that ${\mathcal{CP}}$ is extensible so $z$
has at most $k-1$ parameters, which do not occur in ${\mathcal{CP}}'$. \qed

Take two combinatorial torus cube packings ${\mathcal{CP}}=(z^{i}+[0,1]^{n})_{1\leq i\leq m}$ and ${\mathcal{CP}}'=({z'}^{j}+[0,1]^{n'})_{1\leq j\leq m'}$.
Denote by $({z'}^{i,j}+[0,1]^{n'})_{1\leq j\leq m'}$ with $1\leq i\leq m$
$m$ independent copies of ${\mathcal{CP}}'$; that is every parameter $t'_k$ of
${z'}^j$ is replaced by a parameter $t'_{i,k}$ in ${z'}^{i,j}$.
One defines the combinatorial torus cube packing ${\mathcal{CP}}\ltimes {\mathcal{CP}}'$ by
\begin{equation*}
(z^{i}, {z'}^{i,j})+[0,1]^{n+n'}\mbox{~for~}1\leq i\leq m\mbox{~and~}1\leq j\leq m'.
\end{equation*}
Denote by ${\mathcal{CP}}_1$ the $1$-dimensional combinatorial packing formed by  $(t+[0,1], t+1+[0,1])$.
The combinatorial cube packings ${\mathcal{CP}}_1\ltimes {\mathcal{CP}}_1$ and
${\mathcal{CP}}_1\ltimes ({\mathcal{CP}}_1\ltimes {\mathcal{CP}}_1)$ are
the ones on the left of Figure \ref{2dimContinuousCubePacking}
and \ref{EnumerationDim3}, respectively.
Note that in general ${\mathcal{CP}}\ltimes {\mathcal{CP}}'$ is not
isomorphic to ${\mathcal{CP}}'\ltimes {\mathcal{CP}}$.

\begin{theorem}\label{MostBeautifulTheorem}
Let ${\mathcal{CP}}$ and ${\mathcal{CP}}'$ be two combinatorial torus cube
packings of dimension $n$ and $n'$, respectively.

\begin{itemize}
\item[(i)] $m({\mathcal{CP}}\ltimes {\mathcal{CP}}')=m({\mathcal{CP}}) m({\mathcal{CP}}')$ and $N({\mathcal{CP}}\ltimes {\mathcal{CP}}')=N({\mathcal{CP}})+m({\mathcal{CP}}) N({\mathcal{CP}}')$.

\item[(ii)] ${\mathcal{CP}}\ltimes {\mathcal{CP}}'$ is extensible if and only if ${\mathcal{CP}}$ and ${\mathcal{CP}}'$ are extensible.

\item[(iii)] If ${\mathcal{CP}}$ and ${\mathcal{CP}}'$ are obtained with strictly positive probability and ${\mathcal{CP}}$ is non-extensible then ${\mathcal{CP}}\ltimes {\mathcal{CP}}'$ is attained with strictly positive probability.

\item[(iv)] One has $f^T_{\infty}(n+m)\leq f^T_{\infty}(n)f^T_{\infty}(m)$ and $f^T_{>0,\infty}(n+m)\leq f^T_{>0,\infty}(n) f^T_{>0,\infty}(m)$.
\end{itemize}

\end{theorem}
\proof Denote by $(z^i+[0,1]^{n})_{1\leq i\leq m}$ and by $({z'}^{j}+[0,1]^{n'})_{1\leq j\leq m'}$
the cubes of ${\mathcal{CP}}$ and ${\mathcal{CP}}'$ obtained in this order, i.e.
first $z^1+[0,1]^n$, then $z^2+[0,1]^n$ and so on.
Assertion (i) follows by simple counting.

If ${\mathcal{CP}}$, respectively ${\mathcal{CP}}'$ is extensible
to ${\mathcal{CP}}\cup \{C\}$, ${\mathcal{CP}}'\cup \{C'\}$
then  ${\mathcal{CP}}\ltimes {\mathcal{CP}}'$ is extensible to 
$({\mathcal{CP}}\cup \{C\})\ltimes {\mathcal{CP}}'$, respectively
${\mathcal{CP}}\ltimes ({\mathcal{CP}}'\cup \{C'\})$ and so extensible.
Suppose now that ${\mathcal{CP}}$ and ${\mathcal{CP}}'$ are non-extensible
and take a cube $z+[0,1]^{n+n'}$ with $z$ expressed in terms
of the parameters of ${\mathcal{CP}}\ltimes {\mathcal{CP}}'$.
Then the cube $(z_1,\dots, z_n)+[0,1]^n$ overlaps with one cube
of ${\mathcal{CP}}$, say $z^i+[0,1]^n$.
Also $(z_{n+1}, \dots, z_{n+n'})+[0,1]^{n'}$ overlaps with one cube
of ${\mathcal{CP}}'$, say ${z'}^j+[0,1]^n$.
So, $z+[0,1]^{n+n'}$ overlaps
with the cube $(z^i, {z'}^{i,j})+[0,1]^{n+n'}$ and
${\mathcal{CP}}\ltimes {\mathcal{CP}}'$ is non-extensible, establishing (ii).

A priori there is no simple relation between $p({\mathcal{CP}}\ltimes {\mathcal{CP}}', \infty)$ and
$p({\mathcal{CP}}, \infty)$, $p({\mathcal{CP}}', \infty)$.
But we will prove that if 
$p({\mathcal{CP}}, \infty)>0$, $p({\mathcal{CP}}', \infty)>0$ and
${\mathcal{CP}}$ is not extensible then
$p({\mathcal{CP}}\ltimes {\mathcal{CP}}', \infty)>0$.
That is, to prove (iii) we have to provide one path, among possible many,
in the random sequential cube packing process to obtain
${\mathcal{CP}}\ltimes {\mathcal{CP}}'$ with strictly positive probability
from some corresponding paths of ${\mathcal{CP}}$ and ${\mathcal{CP}}'$.
We first prove that we can obtain the cubes
$((z^{i}, {z'}^{i,1})+[0,1]^{n+n'})_{1\leq i\leq m}$ with strictly
positive probability in this order.
Suppose that we add a cube $z+[0,1]^{n+n'}$ after the cubes $(z^{i'}, {z'}^{i',1})+[0,1]^{n+n'}$ with $i'<i$.
If we choose a coordinate $k\in \{n+1,\dots,n+n'\}$ such that
$z_k=(z^{i'}, {z'}^{i',1})_k+1$ for some $i'<i$ then we still have to
choose a coordinate for all other cubes.
This is because all parameters in $({z'}^{i,1})_{1\leq i\leq m}$ are distinct.
So, we do not gain anything in terms of dimension 
by choosing $k\in \{n+1,\dots,n+n'\}$
and the choice $(z^i, {z'}^{i,1})$ has the same
or higher dimension.
So, we can get the cubes $((z^i, {z'}^{i,1})+[0,1]^{n+n'})_{1\leq i\leq m}$
with strictly positive probability.

Suppose that we have the cubes $(z^i, {z'}^{i,j})+[0,1]^{n+n'}$
for $1\leq i\leq m$ and $1\leq j\leq m'_0$. We will prove by induction
that we can add the cubes $((z^i, {z'}^{i, m'_0+1})+[0,1]^{n+n'})_{1\leq i\leq m}$.
Denote by $n'_{m'_0}\leq n'-1$ the dimension of choices in the 
combinatorial torus cube packing $({z'}^{j}+[0,1]^{n'})_{1\leq j\leq m'_0}$.

Let $z+[0,1]^{n+n'}$ be a cube, which we want to add to the
existing cube packing.
Denote by ${\mathcal S}_z$ the set of $i$ such that $z+[0,1]^{n+n'}$
does not overlap with $(z^i, {z'}^{i,j})+[0,1]^{n+n'}$
on a coordinate $k\leq n$.
The fact that $z+[0,1]^{n+n'}$ does not overlap with the cubes 
$(z^i, {z'}^{i,j})+[0,1]^{n+n'}$ fixes $n'-n'_{m'_0}$ coordinates of $z$.
If $i\not= i'$ then the parameters in ${z'}^{i,j}$ and ${z'}^{i',j'}$ are
different; this means that $(n'-n'_{m'_0})|{\mathcal S}_z|$ components
of $z$ are determined.
Therefore, since ${\mathcal{CP}}$ is non-extensible, we can use 
Lemma \ref{TramLemma}.(ii) and so get the following estimate
on the dimension $D$ of choices:
\begin{equation}\label{OneYearWork}
\begin{array}{rcl}
D  &\leq& \{n' - (n'-n'_{m'_0})|{\mathcal S}_z|\} + \{|{\mathcal S}_z|-1\}\\
   &\leq& n'_{m'_0} - (n'-n'_{m'_0}-1)\{ |{\mathcal S}_z|-1\}\\
   &\leq& n'_{m'_0}.
\end{array}
\end{equation}
We conclude that we cannot do better in terms of dimension
than adding the cubes
$((z^i, {z'}^{i,m'_0+1})+[0,1]^{n+n'})_{1\leq i\leq m}$, which we do.
So we have a path $p$ with $p({\mathcal{CP}}\ltimes {\mathcal{CP}}', p, \infty)>0$ which proves that ${\mathcal{CP}}\ltimes {\mathcal{CP}}'$ is obtained with strictly positive probability.

(iv) follows immediately from (iii) and (ii). \qed

There exist cube packings ${\mathcal{CP}}$, ${\mathcal{CP}}'$ obtained
with strictly positive probability such that
$p({\mathcal{CP}}\ltimes{\mathcal{CP}}', \infty)>0$, which shows that the
hypothesis ${\mathcal{CP}}$ non-extensible is necessary in (iii).

The third $3$-dimensional cube packings of Figure \ref{EnumerationDim3}, 
named rod packing has the cubes $(h^i+[0,1]^3)_{1\leq i\leq 8}$ with the
following $h^i$:
\begin{equation*}
\begin{array}{rcccrccc}
h^1=(&t_1,   &t_2,   &t_3)    &h^5=(&t_6+1, &t_2+1, &t_5+1)\\
h^2=(&t_1+1, &t_4,   &t_5)    &h^6=(&t_1,   &t_2,   &t_3+1)\\
h^3=(&t_6,   &t_2+1, &t_5+1)  &h^7=(&t_1+1, &t_2,   &t_5+1)\\
h^4=(&t_1+1, &t_4+1, &t_5)    &h^8=(&t_1,   &t_2+1, &t_5)\\
\end{array}
\end{equation*}
Taking $8$ $(n-3)$-dimensional combinatorial torus cube-tilings $(w^{i, j})_{1\leq j\leq 2^{n-3}}$ with $1\leq i\leq 8$, one defines a $n$-dimensional {\em rod tiling} combinatorial cube packing
\begin{equation*}
(z^i, w^{i,j})+[0,1]^n\mbox{~~for~~}1\leq i\leq 8\mbox{~~and~~}1\leq j\leq 2^{n-3}.
\end{equation*}

\begin{theorem}
The probability of obtaining a rod tiling is
\begin{equation*}
p^{15}_1\times q_{n-3}^8
\end{equation*}
where $q_n$ is the probability of obtaining a $n$-dimensional cube-tiling and $p^{15}_1$ is a rational function of $n$.
\end{theorem}
\proof Up to equivalence, one can assume that in the random-cube packing process, one puts 
\begin{equation*}
z^1=(h^1, w^{1,1})=(t_1, t_2, t_3, \dots)\mbox{~~and~~}z^2=(h^2, w^{2,1})=(t_1+1, t_4, t_5, \dots).
\end{equation*}
Then there are $n(n-1)$ possible choices for the next cube,
$2(n-1)$ of them are respecting the lamination.
So, there are $(n-2)(n-1)$ choices which do not respect the
lamination and their probability is $p^3_1=\frac{n-2}{n}$.
Without loss of generality, we can assume that one has
\begin{equation*}
z^3=(h^3,w^{3,1})=(t_6, t_2+1, t_5+1, \dots).
\end{equation*}
In the next $5$ stages we add cubes with $n-3$ new parameters each.
We have more than one type to consider under equivalence and we
need to determine the total number of possibilities in order
to compute the probabilities.

For the cube $z^4+[0,1]^n$ we should have three integers $i_1$, $i_2$, $i_3$
such that $z^4_{i_j}\equiv z^j_{i_{j}}+1\pmod 2$. Necessarily, the $i_j$ are
all distinct, which gives $n(n-1)(n-2)$ possibilities.
There are exactly $6$ possibilities with $i_j\leq 3$.
One of them corresponds to the non-extensible cube packing
of Figure \ref{EnumerationDim3} on the first $3$ coordinates which the
$5$ others have a non-zero probability of being extended to the rod tiling.
When computing later probabilities, we used the automorphism group of the
existing configuration and gather the possibilities of extension into
orbits.
At the fourth stage, the $5$ possibilities split into two orbits:
\begin{enumerate}
\item $O^4_1$: $(h^i, w^{i,1})$ for $i\in \{1,2,3,4\}$ with $p^4_1=p^3_1\frac{3}{n(n-1)(n-2)}$,
\item $O^4_2$: $(h^i, w^{i,1})$ for $i\in \{1,2,3,7\}$ with $p^4_2=p^3_1\frac{2}{n(n-1)(n-2)}$; write $\Delta^4_2=3(n-3)(n-4)+3(n-3)+4$ the number of possibilities of adding a cube to the packing $((h^i, w^{i,1})+[0,1]^{n})_{i\in \{1,2,3,7\}}$.
\end{enumerate}
When adding a fifth cube one finds the following cases up to equivalence:
\begin{enumerate}
\item $O^5_1$: $(h^i, w^{i,1})$ for $i\in \{1,2,3,4,5\}$ with $p^5_1=p^4_1\frac{2}{2(n-1)(n-2)}$,
\item $O^5_2$: $(h^i, w^{i,1})$ for $i\in \{1,2,3,4,7\}$ with $p^5_2=p^4_1\frac{2}{2(n-1)(n-2)}+p^4_2\frac{3}{\Delta^4_2}$,
\item $O^5_3$: $(h^i, w^{i,1})$ for $i\in \{1,2,3,7,8\}$ with $p^5_3=p^4_2\frac{1}{\Delta^4_2}$.
\end{enumerate}
When adding a sixth cube one finds the following cases up to equivalence:
\begin{enumerate}
\item $O^6_1$: $(h^i, w^{i,1})$ for $i\in \{1,2,3,4,5,6\}$ with $p^6_1=p^5_1\frac{1}{3(n-2)}$,
\item $O^6_2$: $(h^i, w^{i,1})$ for $i\in \{1,2,3,4,5,7\}$ with $p^6_2=p^5_1\frac{2}{3(n-2)}+p^5_2\frac{2}{n(n-2)}$,
\item $O^6_3$: $(h^i, w^{i,1})$ for $i\in \{1,2,3,4,7,8\}$ with $p^6_3=p^5_2\frac{1}{n(n-2)}+p^5_3\frac{3}{3(n-2)}$.
\end{enumerate}
When adding a seventh cube one finds the following cases up to equivalence:
\begin{enumerate}
\item $O^7_1$: $(h^i, w^{i,1})$ for $i\in \{1,2,3,4,5,6,7\}$ with $p^7_1=p^6_1+p^6_2\frac{1}{n-1}$,
\item $O^7_2$: $(h^i, w^{i,1})$ for $i\in \{1,2,3,4,5,7,8\}$ with $p^7_2=p^6_2\frac{1}{n-1}+p^6_3\frac{2}{2(n-2)}$.
\end{enumerate}
The combinatorial cube packing of eight cubes $((h^i,w^{i,1})+[0,1]^n)_{1\leq i\leq 8}$ is then obtained
with probability $p^8_1=p^7_1+p^7_2\frac{1}{n-2}$.

Then we add cubes in dimension $n-4$ following in fact the construction of
Theorem \ref{MostBeautifulTheorem}. 
The parameters $t_3$, $t_4$ and $t_6$ appear
only two times in the cube packing for the rods, which contain $6$ cubes in
total.
So, when one adds cubes we have $8(n-3)$ choices respecting the cube packing,
i.e. of the form $z^9=(h^i, w^{i, 2})$ with
$w^{i,2}_j\equiv w^{i,1}_j\pmod 2$ for some
$1\leq j\leq n-3$.
We also have $3(n-3)(n-4)$ choices not respecting the rod tiling structure, i.e.
of the form $z^9=(k^i, w)$ with $k^i$ being one of $h^i$ for $1\leq i\leq 3$
with $t_3$, $t_4$ or $t_6$ replaced by another parameter.
But after adding a cube $(h^i, w^{i,2})+[0,1]^{n}$
with $h^i$ containing $t_3$, $t_4$ or $t_6$ this phenomenon cannot occur.
Below a {\em type} $T^h_r$ of probability $p^h_r$ is a packing formed by
the $8$ vectors $(h^i, w^{i,1})_{1\leq i\leq 8}$ and $h-8$ vectors 
of the form $(h^i, w^{i,2})$ amongst which $r$ of the parameters $t_3$, 
$t_4$ or $t_6$ do not occur.
Note that there may be several non-equivalent cube packings with the same
type but this is not important since they have the same numbers of
possibilities.

Adding $9^{th}$ cube one gets:
\begin{enumerate}
\item $T^9_3$, $p^9_1=p^8_1\frac{2(n-3)}{8(n-3)+3(n-3)(n-4)}$,
\item $T^9_2$, $p^9_2=p^8_1\frac{6(n-3)}{8(n-3)+3(n-3)(n-4)}$.
\end{enumerate}
Adding $10^{th}$ cube one gets:
\begin{enumerate}
\item $T^{10}_3$, $p^{10}_1=p^9_1\frac{n-3}{7(n-3)+3(n-3)(n-4)}$,
\item $T^{10}_2$, $p^{10}_2=p^9_1\frac{6(n-3)}{7(n-3)+3(n-3)(n-4)}+p^9_2\frac{3(n-3)}{7(n-3)+2(n-3)(n-4)}$,
\item $T^{10}_1$, $p^{10}_3=p^9_2\frac{4(n-3)}{7(n-3)+2(n-3)(n-4)}$.
\end{enumerate}
Adding $11^{th}$ cube one gets:
\begin{enumerate}
\item $T^{11}_2$, $p^{11}_1=p^{10}_1\frac{6(n-3)}{6(n-3)+3(n-3)(n-4)}+p^{10}_2\frac{2(n-3)}{6(n-3)+2(n-3)(n-4)}$,
\item $T^{11}_1$, $p^{11}_2=p^{10}_2\frac{4(n-3)}{6(n-3)+2(n-3)(n-4)}+p^{10}_3\frac{4(n-3)}{6(n-3)+(n-3)(n-4)}$,
\item $T^{11}_0$, $p^{11}_3=p^{10}_3\frac{2(n-3)}{6(n-3)+2(n-3)(n-4)}$.
\end{enumerate}
Adding $12^{th}$ cube one gets:
\begin{enumerate}
\item $T^{12}_2$, $p^{12}_1=p^{11}_1\frac{n-3}{5(n-3)+2(n-3)(n-4)}$,
\item $T^{12}_1$, $p^{12}_2=p^{11}_1\frac{4(n-3)}{5(n-3)+2(n-3)(n-4)}+p^{11}_2\frac{3(n-3)}{5(n-3)+2(n-3)(n-4)}$,
\item $T^{12}_0$, $p^{12}_3=p^{11}_2\frac{2(n-3)}{5(n-3)+2(n-3)(n-4)}+p^{11}_3\frac{5(n-3)}{5(n-3)+2(n-3)(n-4)}$.
\end{enumerate}
Adding $13^{th}$ cube one gets:
\begin{enumerate}
\item $T^{13}_1$, $p^{13}_1=p^{12}_1\frac{4(n-3)}{4(n-3)+2(n-3)(n-4)}+p^{12}_2\frac{2(n-3)}{4(n-3)+(n-3)(n-4)}$,
\item $T^{13}_0$, $p^{13}_2=p^{12}_2\frac{2(n-3)}{4(n-3)+(n-3)(n-4)}+p^{12}_3\frac{4(n-3)}{4(n-3)}$.
\end{enumerate}
Adding $14^{th}$ cube one gets:
\begin{enumerate}
\item $T^{14}_1$, $p^{14}_1=p^{13}_1\frac{(n-3)}{3(n-3)+(n-3)(n-4)}$,
\item $T^{14}_0$, $p^{13}_2=p^{13}_1\frac{2(n-3)}{3(n-3)+(n-3)(n-4)}+p^{13}_2\frac{3(n-3)}{3(n-3)}$.
\end{enumerate}
Adding $15^{th}$ cube one gets:
\begin{enumerate}
\item $T^{15}_0$, $p^{15}_1=p^{14}_1\frac{2(n-3)}{2(n-3)+(n-3)(n-4)}+p^{14}_2$.
\end{enumerate}
After that if we add a cube $z+[0,1]^n$, then necessarily $z$ is of the form
$(h^i, w)$. 
So, we have $8$ different $(n-3)$-dimensional cube packing problems
show up and the probability is $p^{15}_1 q_{n-3}^8$. \qed

A combinatorial torus cube packing ${\mathcal{CP}}$ is called {\em laminated} if there exist a
coordinate $j$ and a parameter $t$ such that for every cube $z+[0,1]^n$
of ${\mathcal{CP}}$ we have $z_j\equiv t\pmod 1$.

\begin{theorem}\label{UpperBoundLamination}
(i) The probability of obtaining a laminated combinatorial cube packing
is $\frac{2}{n}$.\\[1mm]
(ii) For any $n\geq 1$, one has $E(M_{\infty}^T(n)) \leq 2^n (1-\frac{2}{n})+\frac{4}{n}E(M_{\infty}^T(n-1))$.\\
(iii) For any $n\geq 3$, $\frac{1}{2^n}E(M_{\infty}^T(n))\leq 1-\frac{2^{n}}{n!} \frac{1}{24}$
\end{theorem}
\proof Up to equivalence, we can assume that after the first two
steps of the process, we have
\begin{equation*}
z^1=(t_1,\dots, t_n)\mbox{~~and~~}z^2=(t_1+1, t_{n+1},\dots, t_{2n-1}).
\end{equation*}
So, we consider lamination on the first coordinate.
We then consider all possible cubes that can be added. Those cubes should
have one coordinate differing by $1$ with other vectors. This makes
$n(n-1)$ possibilities. If a vector respects the lamination on the first
coordinate then its first coordinate should be equal to $t_1$ or $t_1+1$.
This makes $2(n-1)$ possibilities.
So, the probability of having a family of cube respecting a lamination
at the third step is $\frac{2}{n}$.
But one sees easily that in all further steps, the choices breaking the
lamination have a dimension strictly lower than the one respecting the
lamination, so they do not occur and we get (i).

By separating between laminated and non-laminated combinatorial torus
cube packings, bounding the number of cubes of non-laminated
combinatorial torus cube packings by $2^n$ one obtains
\begin{equation*}
E(M_{\infty}^T(n))\leq (1-\frac{2}{n})\times 2^n+\frac{2}{n} ( E(M_{\infty}^T(n-1))+E(M_{\infty}^T(n-1))),
\end{equation*}
which is (ii). (iii) follows by induction starting from $\frac{1}{8}E(M_{\infty}^T(3))=\frac{35}{36}$ (see Table \ref{DataInformation}). \qed

\section{Properties of non-extensible cube packings}\label{NonExtensibilityQuestions}
\begin{theorem}\label{ExtendibilityNM3}
If a combinatorial torus cube packing has at least $2^n-3$ cubes, then it is extensible.
\end{theorem}
\proof Our proof closely follows \cite{cubetiling} but is different from it.
Take ${\mathcal{CP}}'$ a combinatorial torus cube packing with
$2^n-\alpha$ cubes, $\alpha\leq 3$.
Take $N$ such that $Nb({\mathcal{CP}}, N)>0$ and ${\mathcal{CP}}$ a
discrete cube packing with $\phi({\mathcal{CP}})={\mathcal{CP}}'$.
If ${\mathcal{CP}}$ is extensible then ${\mathcal{CP}}'$ is extensible
as well.

We select $\delta\in \RR$ and denote by $I_j$ the interval
$[\delta+\frac{j}{2}, \delta+\frac{j+1}{2}[$ for $0\leq j\leq 3$.
Denote by $n_{j,k}$ the number of cubes, whose $k$-th coordinate modulo $2$
belong to $I_j$.

All cubes of ${\mathcal{CP}}$, whose $k$-th coordinate belongs to
$I_j$, $I_{j+1}$ form after removal of their $k$-th coordinate
a cube packing of dimension $n-1$, which we denote by ${\mathcal{CP}}_{j,k}$.
We write $n_{j,k}+n_{j+1,k}=2^{n-1}-d_{j,k}$ and obtain the equations
\begin{equation*}
d_{0,k}-d_{1,k}+d_{2,k}-d_{3,k}=0\quad {\rm and}\quad \sum_{j=0}^3 d_{j,k}=2\alpha.
\end{equation*}
We can then write the vector $d_{k}=(d_{0, k}, d_{1,k}, d_{2, k}, d_{3, k})$ in the following way:
\begin{equation*}
d_{k}=c_1(1,1,0,0)+c_2(0,1,1,0)+c_3(0,0,1,1)+c_4(1,0,0,1)
\quad {\rm with}\quad \sum_{j=1}^4 c_j=\alpha
\end{equation*}
and $c_i\in \ZZ^{+}$. This implies $d_{j,k}=c_{j}+c_{j+1}\leq \sum c_j=\alpha$.
This means that the $(n-1)$-dimensional cube packing ${\mathcal{CP}}_{j,k}$
has at least
$2^{n-1}-3$ cubes, so by an induction argument, we conclude that
${\mathcal{CP}}_{j,k}$ is extensible.

Suppose now that the $k$-th coordinate of the cubes in ${\mathcal{CP}}$ have
values $0<\delta_1 < \delta_2 < \dots < \delta_M < 2$.
So, the set of points in the complement of ${\mathcal{CP}}$,
whose $k$-th coordinate belongs to the interval
$[\delta_{i}, \delta_{i+1}[$ with $\delta_{M+1}=\delta_1+2$
can be filled by translates
of the parallelepiped $Paral_k(\alpha)=[0,1[^{k-1}\times [0, \alpha[\times [0,1[^{n-k}$.

Note that as $\delta$ varies, the vector
$d_{k}$ varies as well.
Suppose that for some $i$, we have the $k$-th layer
$[\delta_i, \delta_{i+1}[$ being full and $[\delta_{i-1}, \delta_i[$
containing $x$ translates with $x\leq 3$
of the parallelepiped $Paral_k(\delta_{i+1}-\delta_{i})$.
Then if one selects another coordinate $k'$, all parallelepipeds
$Paral_{k'}(\delta'_{i'+1}-\delta'_{i'})$
filling the hole delimited by the parallelepiped $Paral_k(\delta_{i+1}-\delta_{i})$ will have the same position in the $k$-th coordinate.
This means that they will form $x$ cubes and that the cube packing
is extensible.
This argument solves the case $\alpha=1$, because up to symmetry
$d_{k}=(0,1,1,0)$.

If $\alpha=2$, then the case of vector of coordinate $d_{k}$ being 
equal to symmetry to $(0,2,2,0)$ or $(0,1,2,1)$ is also solved
because we have seen
that a full layer implies that we can fill the hole.
We have the remaining case $(1,1,1,1)$. If the hole of this
cube packing cannot be filled, then we have a structure of this form:
\begin{center}
\begin{picture}(0,0)%
\includegraphics{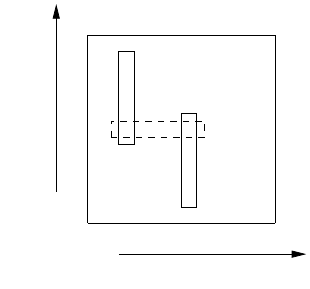}%
\end{picture}%
\setlength{\unitlength}{1973sp}%
\begingroup\makeatletter\ifx\SetFigFont\undefined%
\gdef\SetFigFont#1#2#3#4#5{%
  \reset@font\fontsize{#1}{#2pt}%
  \fontfamily{#3}\fontseries{#4}\fontshape{#5}%
  \selectfont}%
\fi\endgroup%
\begin{picture}(2952,2737)(361,-2186)
\put(376,-436){\makebox(0,0)[lb]{\smash{{\SetFigFont{9}{10.8}{\rmdefault}{\mddefault}{\updefault}{\color[rgb]{0,0,0}$x_{k'}$}%
}}}}
\put(2101,-2086){\makebox(0,0)[lb]{\smash{{\SetFigFont{9}{10.8}{\rmdefault}{\mddefault}{\updefault}{\color[rgb]{0,0,0}$x_{k}$}%
}}}}
\end{picture}%

\end{center}
Selecting another coordinate $k'$, we get that the two parallelepipeds
$z+Paral_k(\delta_{i+1}-\delta_i)$ and $z'+Paral_k(\delta_{i'+1}-\delta_{i'})$
have $z_{l}=z'_{l}$ for $l\not= k, k'$. This is impossible if $n\geq 4$.
So, if $d_{k}=(1,1,1,1)$ for some $k$ and $\delta$, then the hole
can be filled.

If $\alpha=3$, and $d_{k}$, up to symmetry, is equal to $(0,3,3,0)$ or $(0,2,3,1)$ then we have 
a full layer and so we can fill the hole.
If the vector $d_{k}=(2,1,1,2)$
occurs, then by the same argument as for $(1,1,1,1)$ we can fill
the hole. \qed

\begin{proposition}\label{FirstStupidProposition}
(i) Non extensible combinatorial torus cube packings of dimension $n$
have at least $n+1$ cubes.

(ii) If ${\mathcal{CP}}$ is a combinatorial torus cube tiling, 
then in a coordinate $j$ a parameter $t$
occur the same number of times as $t$ and $t+1$.
\end{proposition}
\proof (i) Suppose that a combinatorial torus cube packing ${\mathcal{CP}}$
has $m\leq n$ cubes $(z^i+[0,1]^n)_{1\leq i\leq m}$.
By fixing $z_i=z^i_i+1$ for $i=1, 2, \dots, m$ we get that
the cube $z+[0,1]^{n}$ does not overlap with ${\mathcal{CP}}$.

(ii) Without loss of generality, we can assume that a given parameter $t$
occurs only in one coordinate $k$ as $t$ and $t+1$.
The cubes occurring in the layer $[t,t+1]$, $[t+1,t+2]$ on $j$-th
coordinate are the ones with $x_j=t$, $t+1$; we denote by $V_{t}$ and $V_{t+1}$
their volume.
Now if we interchange $t$ and $t+1$ we still obtain a tiling, so $V_{t}\leq V_{t+1}$ and $V_{t+1}\leq V_t$. So, $V_{t}=V_{t+1}$ and the number of cubes with
$x_j=t$ is equal to the number of cubes with $x_{j}=t+1$. \qed

Take a combinatorial torus cube packing ${\mathcal{CP}}$ obtained with strictly
positive probability.
Let us choose a path $p$ to obtain ${\mathcal{CP}}$.
Denote by $N_{k, p}({\mathcal{CP}})$ the number of cubes obtained
with $k$ new parameters along the path $p$.

\begin{proposition}\label{NumberParameters}
Let ${\mathcal{CP}}$ be a non-extensible combinatorial torus cube packing, 
$p$ a path with $p({\mathcal{CP}}, p, \infty)>0$.

(i) $N_{n, p}({\mathcal{CP}})=1$ and $N_{n-1, p}({\mathcal{CP}})=1$.

(ii) $N_{n-2, p}({\mathcal{CP}})\leq 2$ and 
$N_{n-2, p}({\mathcal{CP}}) = 2$ if and only if ${\mathcal{CP}}$ is laminated.

(iii) One has $N_{k,p}({\mathcal{CP}})\geq 1$ for $0\leq k\leq n$.

(iv) $N({\mathcal{CP}})=\sum_{k=0}^{n} kN_{k, p}({\mathcal{CP}})\geq \frac{n(n+1)}{2}$.

(v) If $N({\mathcal{CP}})=\frac{n(n+1)}{2}$ then $N_{k,p}({\mathcal{CP}})=1$ for $k\geq 1$.
\end{proposition}
\proof The first cube $z^1+[0,1]^n$ has $n$ new parameter, but the
second cube $z^2+[0,1]^n$ should not overlap with the first one
so it has $n-1$ parameters and $N_{n,p}({\mathcal{CP}})=1$.
Without loss of generality, we can assume that $z^1=(t_1, \dots, t_n)$
and $z^2=(t_1+1, t_{n+1}, \dots, t_{2n-1})$.
When adding the third cube $z^3+[0,1]^n$, we have to set up $2$
coordinates depending on the parameters $t_i$, $i\leq 2n-1$ thus
$N_{n-1, p}({\mathcal{CP}})=1$.

If $z^3_1=t_1$ or $t_1+1$ then we have a laminated cube packing, we can
add a cube with $n-2$ parameters and $N_{n-2, p}({\mathcal{CP}})=2$.
Otherwise, we do not have a laminated cube packing, three coordinates
of $z^4$ need
to be expressed in terms of preceding cubes and thus
$N_{n-2, p}({\mathcal{CP}})=1$.

(iii) The proof is by induction;
suppose one has put $m'=\sum_{l=k}^{n} N_{l, p}({\mathcal{CP}})$ cubes.
Then the cube $z^{m'}+[0,1]^n$ has $k$ new parameters $t'_1$, \dots, $t'_k$
in coordinates $i_1$, \dots, $i_k$.
The cube $C=z+[0,1]^n$ with $z_{i_1}=t'_1+1$ and 
$z_i=z^{m'}_i$ for $i\notin \{i_1, \dots, i_k\}$
has $k-1$ free coordinates $\{i_2, \dots,i_k\}$
and thus $k-1$ new parameters.
So, $N_{k-1, p}({\mathcal{CP}})\geq 1$.

(iv) and (v) are elementary. \qed

\begin{conjecture}\label{Power2conjecture}
Let ${\mathcal{CP}}$ be a combinatorial torus cube packing and $p$ a path with $p({\mathcal{CP}}, p, \infty)>0$.

(i) For all $k\geq 1$ one has $\sum_{k=0}^l N_{n-k, p}({\mathcal{CP}})\leq 2^l$.

(ii) $N({\mathcal{CP}})\leq 2^n-1$; 
if $N({\mathcal{CP}}) = 2^n-1$, then ${\mathcal{CP}}$ is obtained
via a lamination construction.
\end{conjecture}

A {\em perfect matching} of a graph $G$ is a set ${\mathcal M}$ of edges
such that every vertex of $G$ belongs to exactly one edge of ${\mathcal M}$.
A {\em $1$-factorization} of a graph $G$ is a set of perfect matchings,
which partitions the edge set of $G$.
The graph $\mathsf{K}_4$ has one $1$-factorization; 
the graph $\mathsf{K}_6$ has, up to isomorphism, exactly
one $1$-factorization with symmetry group $\Sym(5)$.

\begin{proposition}\label{MinimalNrParam}
Let ${\mathcal{CP}}$ be a non-extensible combinatorial torus cube packing.

(i) If $n$ is even then ${\mathcal{CP}}$ has at least $n+2$ cubes.

(ii) If $n$ is odd and ${\mathcal{CP}}$ has $n+1$ cubes then
$N({\mathcal{CP}})=\frac{n(n+1)}{2}$. Fix a coordinate $j$ and a parameter
$t$ occurring in at least one cube.
Then the number of cubes containing $t$, respectively $t+1$
in coordinate $j$ is exactly $1$.

(iii) If $n$ is odd then isomorphism classes of non-extensible
combinatorial torus cube packings with $n+1$ cubes are in one
to one correspondence
with isomorphism classes of $1$-factorizations of $\mathsf{K}_{n+1}$.

(iv) If $n$ is odd then the non-extensible combinatorial torus
cube packings with $n+1$ cubes
are obtained with strictly positive probability and
$f^T_{\infty}(n)=f^T_{>0,\infty}(n)=n+1$.

\end{proposition}
\proof We take a non-extensible cube packing ${\mathcal{CP}}$ with $n+1$ cubes.
Suppose that for a coordinate $j$ we have two cubes $z^i+[0,1]^n$ and
$z^{i'}+[0,1]^n$ with $z^i_{j}=z^{i'}_j=t$.
If a vector $z$ has $z_{j}=t+1$, then $z+[0,1]^n$ does not overlap with
$z^i+[0,1]^n$ and $z^{i'}+[0,1]^n$. There are $n-1$ remaining cubes to
which $z+[0,1]^n$ should not overlap but we have $n-1$ remaining coordinates
so it is possible to choose the coordinates of $z$ so that $z+[0,1]^n$ does
not overlap with ${\mathcal{CP}}$.
This is impossible, therefore parameters appear always
at most $1$ time as $t$ and at most one
time as $t+1$ in a given coordinate.

By Lemma \ref{TramLemma} every parameter $t$ appear also as $t+1$.
So, every parameter $t$ appears one time as $t$ and one time as $t+1$.
This implies that we have an even number of cubes and so (i).
Every coordinate has $\frac{n+1}{2}$ parameters,
which gives $\frac{n(n+1)}{2}$ parameters and so (ii).

(iii) Assertion (ii) implies that any two cubes $C_i$ and $C_{i'}$ of
${\mathcal{CP}}$ have exactly one
coordinate on which they differ by $1$. So, every coordinate correspond
to a perfect matching and the set of $n$ coordinates to the $1$-factorization.

(iv) Since parameters $t$ appear only one time as $t$ and $t+1$, the
dimension of choices after $k$ cubes are put is $n-k$
and one sees that such a cube packing is obtained with
strictly positive probability.
The existence of $1$-factorization of $\mathsf{K}_{2p}$
(see, for example, \cite{walecki,harary})
gives $f^T_{\infty}(n)\leq f^T_{>0,\infty}(n)\leq n+1$.
Combined with Theorem \ref{FirstStupidProposition}.i, 
we have the result. \qed

\begin{conjecture}
If $n$ is even then there exist non-extensible combinatorial torus
cube packings with $n+2$ cubes and $\frac{n(n+1)}{2}$ parameters.
\end{conjecture}

In dimension $4$ there is a unique cube packing (obtained with
probability $\frac{1}{480}$) satisfying this conjecture:
\begin{equation*}
\left(\begin{array}{cccc}
t_{1}   & t_{2}   & t_{3}   & t_{4}  \\
t_{5}   & t_{6}   & t_{7}   & t_{4}+1\\
t_{1}+1 & t_{8}   & t_{7}+1 & t_{9}  \\
t_{5}+1 & t_{8}+1 & t_{3}+1 & t_{10}  \\
t_{1}+1 & t_{6}+1 & t_{7}   & t_{10}+1\\
t_{5}   & t_{2}+1 & t_{7}+1 & t_{9}+1
\end{array}\right)
\end{equation*}

\begin{figure}
\newenvironment{MYmatrix}{
\setlength{\arraycolsep}{1pt}
\begin{matrix}}{\end{matrix}}

\begin{center}
\input{flte8_6.tex}
\end{center}
\caption{The non-extensible $6$-dimensional combinatorial cube packings with $8$ cubes and at least $21$ parameters}
\label{The9structures}
\end{figure}

\begin{proposition}\label{6dimCases}
(i) There are $9$ isomorphism types of non-extensible combinatorial torus
cube packings in dimension $6$ with $8$ cubes and
at least $21$ parameters (see Figure \ref{The9structures});
they are not obtained with strictly positive probability.

(ii) $8=f^T_{\infty}(6) < f^T_{>0,\infty}(6)$.
\end{proposition}
\proof (ii) follows immediately from (i).
The enumeration problem in (i) is solved in the following way:
instead of adding cube after cube like in the random cube packing process,
we add coordinate after coordinate in all possible ways and reduce
by isomorphism. The computation returns the listed combinatorial torus
cube packings.
Given a combinatorial cube packing ${\mathcal{CP}}$ in order to prove that
$p({\mathcal{CP}}, \infty)=0$, we consider all ($8!$) possible paths $p$
and see that for all of them $p({\mathcal{CP}}, p, \infty)=0$. \qed

\begin{proposition}\label{StrangeDiscreteCases}
If $n=3$, $5$, $7$, $9$, then there exist a combinatorial torus cube tiling
obtained with strictly positive probability and $\frac{n(n+1)}{2}$ parameters.
\end{proposition}
\proof If $n$ is odd consider the matrix $H_n=m_{i,j}$ with all
elements satisfying $m_{i+k, i}=m_{i-k, i}+1$ for
$1\leq k\leq\frac{n-1}{2}$, the addition being modulo $n$.
The matrix for $n=5$ is
\begin{equation*}
H_5=\left(\begin{array}{ccccc}
t_1      & t_7+1    & t_{13}+1 & t_{14}   & t_{10}\\
t_6      & t_2      & t_8+1    & t_{14}+1 & t_{15}\\
t_{11}   & t_7      & t_3      & t_9+1    & t_{15}+1\\
t_{11}+1 & t_{12}   & t_8      & t_4      & t_{10}+1\\
t_6+1    & t_{12}+1 & t_{13}   & t_9      & t_5\\
\end{array}\right).
\end{equation*}
Then form the combinatorial cube packing with the cubes $(z^i+[0,1]^n)_{1\leq i\leq n}$
and $z^i$ being the $i$-th row of $H_n$.
It is easy to see that the number of parameters of cubes,
which we can add after $z^i$ is $n-i$.
So, those first $n$ cubes are attained with the minimal
number $\frac{n(n+1)}{2}$ of parameters and with strictly
positive probability.
If a cube $z+[0,1]^n$ is non-overlapping with $z^i+[0,1]^n$ for $i\leq n$
then there exist $\sigma(i)\in \{1,\dots,n\}$
such that $z_{\sigma(i)}=z^i_{\sigma(i)}+1$.
If $i\not= i'$ then $\sigma(i)\not= \sigma(i')$, which proves that
$\sigma\in \Sym(n)$. We also add the $n$ cubes corresponding to the matrix
$H_n+Id_n$.
So, there are $n!$ possibilities for adding new cubes and we
need to prove that we can select $2^n-2n$ non-overlapping cubes
amongst them.

The symmetry group of the $n$ cubes $(z^i+[0,1]^n)_{1\leq i\leq n}$
is the dihedral group $D_{2n}$ with $2n$ elements.
It acts on $\Sym(n)$ by conjugation and so we simply
need to list the relevant set of inequivalent permutations in order
to describe the corresponding cube packings.
See Table \ref{TheComputerComputation} for the found permutation
for $n=3$, $5$, $7$, $9$. \qed

The cube packing of above theorem was obtained for $n=5$
by random method, i.e., adding cube whenever possible by
choosing at random.
Then the packings for $n=7$ and $9$ were built using the matrix $H_n$
and consideration of all possibilities invariant under the dihedral
group $D_{2n}$ by computer. But for $n=11$ this method does not work.
It would be interesting to know in which dimensions $n$ 
combinatorial torus cube tilings
with $\frac{n(n+1)}{2}$ parameters do exist.

\begin{table}
\begin{center}
\begin{tabular}{|c|ccc|}
\hline
$n=3$ & $(1,2,3)$\\
\hline
$n=5$ & $(1,2,3,4,5)$ & $(1,2)(3,5,4)$ & $(1,4,5,3,2)$\\
\hline
$n=7$ & $(1,2,3,4,5,6,7)$ & $(1,7)(2,5,4,3,6)$ & $(1,6,2,5,4,3,7)$\\
      & $(1,7)(2,5,6)(3,4)$ & $(1,6,2,3,7)(4,5)$ & $(1,7)(2,3,4,5,6)$\\
      & $(1,3,7)(2,6)(4,5)$ & $(1,3,7)(2,5,4,6)$ & $(1,5,4,3,2,6,7)$\\
\hline
$n=9$ &$(1,2,3,4,5,6,7,8,9)$ & $(1,6,7,4,3,5,9)(2,8)$ & $(1,5,6,7,4,3,9)(2,8)$\\
      &$(1,5,9)(2,8)(3,6,7,4)$ & $(1,9)(2,5,4,3,6,7,8)$ & $(1,6,7,8,2,5,4,3,9)$\\
      &$(1,5,4,3,9)(2,8)(6,7)$ & $(1,9)(2,5,6,7,8)(3,4)$ & $(1,9)(2,3,6,5,4,7,8)$\\
      &$(1,6,5,4,7,8,2,3,9)$ & $(1,9)(2,3,6,7,8)(4,5)$ & $(1,6,7,8,2,3,9)(4,5)$\\
      &$(1,9)(2,3,4,7,6,5,8)$ & $(1,9)(2,3,4,5,6,7,8)$ & $(1,9)(2,3,6)(4,7,8,5)$\\
      &$(1,6,2,3,9)(4,7,8,5)$ & $(1,5,4,7,8,6,2,3,9)$ & $(1,8,3,7,6,4,9)(2,5)$\\
      &$(1,7,6,4,9)(2,5)(3,8)$ & $(1,7,6,9)(2,8,3,4,5)$ & $(1,4,5,2,8,3,7,6,9)$\\
      &$(1,4,9)(2,5)(3,7,6,8)$ & $(1,4,8,3,7,6,9)(2,5)$ & $(1,7,6,3,4,5,2,8,9)$\\
      &$(1,4,5,2,8,9)(3,7,6)$ & $(1,4,9)(2,5)(3,8,7,6)$ & $(1,3,7,6,9)(2,8,4,5)$\\
      &$(1,7,6,5,4,3,2,8,9)$ & $(1,9)(2,5,8)(3,4)(6,7)$\\
\hline
\end{tabular}
\end{center}
\caption{List of permutation describing combinatorial torus cube tilings with $\frac{n(n+1)}{2}$ parameters in dimension $3$, $5$, $7$, $9$}
\label{TheComputerComputation}
\end{table}

\section{Acknowledgments}
We thank Luis Goddyn and the anonymous referees for helpful comments.

\end{document}